\documentstyle[12pt]{article}

\newcommand{\app}[1]{\setcounter{section}{0}
\setcounter{equation}{0} \renewcommand{\thesection}
{\Alph{section}}\section{#1}}
\newcommand{\sect}[1]{\setcounter{equation}{0}\section{#1}}
\newcommand{\subsect}[1]{\setcounter{equation}{0}\subsection{#1}}
\def\noi{\noindent}
\def\ar{\begin{array}{rcl}}
\def\an{\end{array}}

\newcommand{\eq}{\begin{equation}}
\newcommand{\eqa}{\begin{eqnarray}}
\newcommand{\en}{\end{equation}}
\newcommand{\ena}{\end{eqnarray}}
\def\1{{\bf 1}}

\def\F{\mbox{$\cal F$\,}}

\def\Fp{{\cal F}_{21}}

\def\ot{\otimes}
\def\id{\mbox{id}}
\def\ie{\mbox{\it i.e.\/ }}
\def\eg{\mbox{\it e.g.\/ }}
\def\C{{\cal C}}
\def\up{\uparrow}
\def\idA{{\bf 1}_{\cal A}}
\def\idAh{{\bf 1}_{{\cal A}_h}}

\def\g{\mbox{\bf g\,}}
\def\uqg{\mbox{$U_h{\/\mbox{\bf g}}$ }}

\newcommand{\tr}{\triangleright_h\,}
\newcommand{\trl}{\tilde\triangleright_h\,}
\newcommand{\trc}{\triangleright}

\def\R{\mbox{$\cal R$\,}}

\def\A{\mbox{${\cal A}_{\pm,{\footnotesize\mbox{\bf g\,}},\rho}$}}
\def\Aq{\mbox{${\cal A}^h_{\pm,{\footnotesize\mbox{\bf g\,}},\rho}$}}
\def\Ai{\mbox{${\cal A}^{inv}_{\pm,{\footnotesize\mbox{\bf
g\,}},\rho}$}}
\def\Aiq{\mbox{${\cal A}^{h,inv}_{\pm,{\footnotesize\mbox{\bf
g\,}},\rho}$}}
\def\aa{\mbox{${\sf a}$}}

\def\P{\mbox{$P^F$}}
\def\F{\mbox{$\cal F$}}

\def\z{\hspace*{9mm}}
\def\x{\hspace{3mm}}

\newcommand{\cn}{{\bf C}}
\newcommand{\rn}{{\bf R}}
\newcommand{\zn}{{\bf Z}}
\newcommand{\nn}{{\bf N}}
\newtheorem{prop}{Proposition}
\newtheorem{lemma}{Lemma}
\newtheorem{theorem}{Theorem}

\begin{document}
%\pagestyle{fancy}

%\cfoot{}
%\rhead[\fancyplain{}]{{\thepage}}
                
\begin{titlepage}
\begin{center}

~

%July 1997         \hfill       LMU-TPW 97-07,  ~LBNL-40274\\
\vskip.6in

{\Large \bf Drinfel'd Twist and q-Deforming Maps for\\
            Lie Group Covariant Heisenberg Algebrae}

\vskip.4in

Gaetano Fiore* 

\vskip.25in

{\em Sektion Physik der Ludwig-Maximilians-Universit\"at
M\"unchen\\
Theoretische Physik --- Lehrstuhl Professor Wess\\
Theresienstra\ss e 37, 80333 M\"unchen\\
Federal Republic of Germany}
\end{center}
\vskip1in
\begin{abstract}
Any deformation of a Weyl or Clifford algebra can  be realized 
through a change of generators in the undeformed algebra.
$q$-Deformations 
of Weyl or Clifford algebrae that were covariant under the action of 
a simple Lie algebra $\g$ are characterized by their being 
covariant under the action of the quantum group $\uqg$, $q:=e^h$.
We present a systematic procedure for determining 
all possible corresponding changes of generators, together with
the corresponding realizations of the $\uqg$-action.
The intriguing relation between $\g$-invariants and $\uqg$-invariants
suggests that these changes of generators might be employed
to simplify the dynamics of some $\g$-covariant quantum physical
systems.
\end{abstract}
\vfill
\begin{center}
~
%PACS: ~ 02.20.-a,~ 03.65.Fd,~ 05.30.-d,~ 11.30.-j
\end{center}
\noi \hrule
\vskip.2cm
\noi{\footnotesize *EU-fellow, TMR grant ERBFMBICT960921.
\qquad {\it e-mail: }Gaetano.Fiore@physik.uni-muenchen.de}
\end{titlepage}
\newpage
\setcounter{page}{1}

\renewcommand{\theequation}{\thesection.\arabic{equation}}
\sect{Introduction}

Weyl and Clifford algebrae (respectively denoted by 
${\cal A}_+,{\cal A}_-$ 
in the sequel, and collectively as ``Heisenberg
algebrae'') are at the hearth of quantum physics.
One may ask  if deforming them within the category of associative 
algebrae (i.e. deforming their 
defining commutation relations) yields new physics, 
or at least may be useful to better describe
certain systems in conventional quantum physics.
This question can be divided into an algebraic and a 
representation-theoretic subquestion.  

The first was addressed in the fundamental paper \cite{flato}.
Essentially, it reads: 
is there a formal realization of the elements of the deformed
algebra in terms of elements of the undeformed algebra? 
The answer is affirmative, 
but in general the realization is not explicitly known.
A general result \cite{ducloux} 
regarding the Hochschild cohomology
of the universal enveloping algebra associated to a 
nilpotent Lie group 
states in particular that the first and second 
cohomology groups of any 
Weyl algebra ${\cal A}_+$ 
are trivial. This implies \cite{gerst} that any 
deformation ${\cal A}_+^h$
($h$ denoting the deformation parameter) of the latter
is trivial, in the sense that there exists
an isomorphism of topological algebrae over $\cn[[h]]$
(a ``deforming map'', in the terminology of Ref. \cite{zac}),
$f:{\cal A}_+^h\rightarrow{\cal A}_+[[h]]$, reducing
to the identity in the limit $h=0$
(a concise and effective presentation of these 
results can be found in Sect.'s 1,2 of 
Ref. \cite{mathias}).
Practically this means  that the
generators $\tilde A^i,\tilde A_i^+$ 
of ${\cal A}_+^h$ are mapped by $f$
into power series in $h$ with coefficients in ${\cal A}_+$,
$A^i:=f(\tilde A^i)$, $A^+_i=f(\tilde A_i^+)$, fulfilling
the same deformed commutation relations and
going to the corresponding
generators $a^i,a^+_i$ of ${\cal A}_+$ in the limit $h=0$.

The second subquestion is: do 
also the deformed and undeformed representation theories
coincide? One can already see in some simple model that the answer is
negative, but in the general case, up to our knowledge, 
the relation between the two is an open question.

Given any automorphism 
$g:{\cal A}_+[[h]]\rightarrow{\cal A}_+[[h]]$, $g=\id+O(h)$,
then $g\circ f$ is a new deforming map; conversely, given
two deforming maps $f,f'$, the map $f'\circ f^{-1}$ is an
algebra automorphism. Now, 
by the vanishing of the first cohomology group of ${\cal A}_+$,
all automorphisms of ${\cal A}_+[[h]]$ are `inner',
\ie of the form $g(a)=\alpha\, a \,\alpha^{-1}$. Hence, all
deforming maps can be obtained from one through the
formula
\eq
f_{\alpha}(a):=\alpha f(a)\alpha^{-1} \z\z \alpha=1+O(h)
\in{\cal A}_+^h.
\label{innerauto}
\en

These results apply \cite{mathias} in particular to so-called
``$q$-deformations'' ($q:=e^h$)  of Weyl algebrae which are
covariant under the action of some simple Lie algebra $\g$;
such deformations \cite{puwo,wezu,cawa} 
are matched to the deformation of $U\g$ into the quantum group 
$\uqg$ \cite{dr2},
in the sense that for all $q$ the deformed algebrae are in fact
\uqg-module algebrae\footnote{They should not be confused with the
celebrated Biedenharn-Macfarlane-Hayashi $q$-oscillator 
(super)algebrae \cite{mac}, whose generators 
 are {\it not} 
\uqg-covariant (in spite of the fact that they are usually used
to construct a generalized Jordan-Schwinger realization of \uqg).}: 
the commutation relations among $\tilde A^i,\tilde A_i^+$
are compatible with the quantum group action 
$\trl:\uqg\times{\cal A}_+^h\rightarrow 
{\cal A}_+^h$ as the commutation relations among $a^i,a^+_i$ were
compatible with the classical action 
$\trc:U\g\times {\cal A}_+
\rightarrow {\cal A}_+$.

Also quantum groups admit algebra isomorphisms
$\varphi_h:\uqg\rightarrow U\g[[h]] ~\cite{dr3}$\footnote{
The existence of the latter and their being defined up to
inner automorphisms of $U\g[[h]]$ again
is a consequence of the triviality of the first and second
Hochschild cohomology groups of $U\g$.}.
In spite of the existence of the algebra
isomorphisms  $f,\varphi_h$,
the $\uqg$-module algebra structure
$(\uqg,{\cal A}_+^h,\trl)$ is 
a {\it non-trivial} deformation of the one
$(U\g,{\cal A}_+,\trc)$, \ie 
for no $\varphi_h,f$ the equality
$f\circ \trl=\trc\circ(\varphi_h\times f)$ holds\footnote{
By $\trc$ we mean here actually its linear extension to
$\trc:U\g[[h]]\times {\cal A}_{+,\footnotesize\mbox{\bf g\,},\rho}[[h]]
\rightarrow {\cal A}_{+,\footnotesize\mbox{\bf g\,},\rho}[[h]]$, where
both the domain and codomain have to be understood 
completed in the $h$-adic topology.}.
This is because \uqg itself as a {\it Hopf} algebra
is a {\it non-trivial} deformation of $U\g$, in other words 
all $\varphi_h$'s are algebra but
not coalgebra (and therefore not Hopf algebra) isomorphisms
(this is related to the non-triviality of
the Gerstenhaber-Schack cohomology \cite{gescha}).

Given $f$ we define $\tr$ as the map making 
the following diagram commutative:
\eq
\begin{array}{ccccc}
 \uqg & \times & {\cal A}_{\pm}^h & \, 
 \stackrel{\trl}{-\!\!\!-\!\!\!-\!\!\!- \!\!\!-\!\!\!-
\!\!\!\longrightarrow}\, & {\cal A}_{\pm}^h \cr
 &\updownarrow& \id \times  f  & &\updownarrow\, f \cr
 \uqg & \times & {\cal A}_{\pm}[[h]] & \, 
 \stackrel{\tr}{-\,-\,-\rightarrow}\, & {\cal A}_{\pm}[[h]]\cr
\end{array}
\label{diagram}
\en
(in other words
$\tr:= f\circ\trl\circ(\id\ot f^{-1})$ will realize $\trl$ on
$\A[[h]]$).
In this work we give a systematic procedure to construct all pairs
$(f,\tr)$. It is based on
the properties of the ``Drinfel'd twist'' \cite{dr3}. We first
construct $\tr$, then $f$.

One particular $\tr$ can be naturally constructed 
in a `adjoint-like'  way (Sect. \ref{qcov}). 
To determine a corresponding $f$ it is sufficient to identify in
$\A[[h]]$ appropriate realizations $A^i,A^+_i$  of
$\tilde A^i,\tilde A_i^+$. With this aim in mind,  
we first show [formula (\ref{def3})]
how to construct candidates $A^i,A^+_i$ having
the same transformation properties 
under $\tr$ as the generators $\tilde A^i,\tilde A_i^+$
of $\Aq$ under $\trl$; 
they are determined up to multiplication by some
\g-invariants. By requiring that the $A^i,A_i^+$
also have the same commutation rules as the
$\tilde A^i,\tilde A_i^+$ (Sect. \ref{reqcr}) we find constraints
(\ref{figata1}-\ref{figata3}) on these \g-invariants.
Requiring that the
$*$-structures of ${\cal A}_{\pm}$ realize the $*$-structures of 
${\cal A}_{\pm}^h$ yields further  constraints (Sect. \ref{star}).
The subalgebrae ${\cal A}_{\pm}^{inv}[[h]],{\cal A}_{\pm}^{h,inv}[[h]]$
of $\A[[h]]$ that are invariant respectively under $\trc,\tr$ coincide
(Sect. \ref{inva}), but invariants in the form of polynomials in 
$A^i,A^+_j$ are highly non-polynomial 
(analytic) functions of the classical invariants
and hence of $a^i, a^+_j$, and conversely; we explicitly find 
these functions.
Solving the system of equations (\ref{figata1}-\ref{figata3})
depends on the particular \g~ and on the particular
$U\g$-representation $\rho$ to which the generators $a^+_i$
belong (the $a^i$ necessarily belong to the contragradient
$\rho^\vee$ of $\rho$). We shall denote by 
$\A$ the corresponding Heisenberg algebra and by $\Aq$
its $q$-deformed version. We solve (Sect. \ref{reqcr})
the above equations
for  the well-known  ${\cal A}^h_{\pm,sl(N),\rho_d}$,
${\cal A}^h_{+,so(N),\rho_d}$,
($\rho_d$ will denote the  defining representations of either \g). 
Finally, in Sect. \ref{conclu} we extend the previous
results to all other isomorphisms (\ref{innerauto})
while giving an outlook of the whole construction, we make
some remarks on the representation theory,
and we draw the conclusions. This serves also to better 
clarify the motivations for the present work.

In Ref. \cite{fio} we started the program just
sketched by sticking to the cases of 
triangular deformations of the Hopf algebra $U\g$
(in the present
setting it is recovered by postulating a trivial coassociator)
and of the $q$-deformation ${\cal A}^h_{\pm,sl(2),\rho_d}$.

Examples of $q$-deforming maps  for
Heisenberg algebrae were explicitly 
constructed ``by hand'' in past works
\cite{zachos3,oleg,mathias,wess}\footnote{For multidimensional
Heisenberg algebrae the non-trivial task there was to find
an intermediate transformation 
to a set of mutually commuting {\it pairs}
$\{\alpha^i,\alpha^+_i\}$ of deformed generators; 
these generators $\alpha^i,\alpha^+_i$
are covariant neither under the $U\g$ nor the $U_h\g$ action.}.
Of course these deforming maps are
related to ours by some automorphism (\ref{innerauto});
in Sect. \ref{better} we determine the latter 
for ${\cal A}^h_{+,sl(2),\rho_d}$
and the deforming map found in Ref. \cite{oleg}.
 We underline that
our construction is based instead on universal objects
characterizing the quantum symmetry algebra. This allows
the application of our method \eg to the physically relevant
case that $\rho$ be the direct sum of an arbitrary number of
copies of $\rho_d$.

\renewcommand{\theequation}{\thesection.\arabic{equation}}
\renewcommand{\theequation}{\thesubsection.\arabic{equation}}
\sect{Preliminaries and notation}

Some general remarks before starting.
Although we will always denote 
the generators of the Heisenberg algebrae by 
$a^i,a_i^+,A^i,A_i^+,...$, only the choice of a 
$*$-structure may give them  the meaning of 
creators/annihilators, or coordinates/derivatives, etc. 
(see \eg Sect. \ref{star}).
Given an algebra $B$, $B[[h]]$ will denote 
 the algebra of formal power series in $h\in \cn$
with coefficients belonging to finite-dimensional subspaces
of $B$. Both $B[[h]]$ and tensor 
products like $B[[h]]\ot B[[h]]$ will be understood to be
completeted in the $h$-adic topology. The symbol $U_h\g$ \cite{dr2} 
will denote the algebra on the ring $\cn[[h]]$
underlying the quantum group (completed in the $h$-adic topology).

\subsection{Twisting groups into quantum groups}

Let $H=(U\g,m,\Delta, \varepsilon,S )$ be the 
cocommutative Hopf algebra associated to the universal
enveloping (UE) algebra $U\/ \g$ of a Lie algebra \g. 
The symbol $m$ denotes the multiplication (in the 
sequel it will be dropped in the obvious way 
$m(a\ot b)\equiv ab$, unless explicitly required), whereas
$\Delta, \varepsilon,S $  the comultiplication, counit 
and antipode respectively. Assume that 
$H_h=(\uqg,m_h,\Delta_h,\varepsilon_h,S_h,\R)$ is a quasitriangular
non-cocommutative deformation of $H$ ~\cite{dr2}; 
$h\in\cn$ denotes the deformation parameter, 
$m_h,\Delta_h,\varepsilon_h,S_h$ the deformed multiplication,
comultiplication, counit and antipode respectively, and 
$\R\in\uqg\otimes\uqg$ 
the universal $R$-matrix.

A well-known theorem by Drinfel'd, Proposition 3.16
in Ref. \cite{dr3} (whose results 
are to a certain extent already implicit in preceding works by Kohno 
\cite{ko}), proves, for any simple finite-dimensional Lie algebra \g,
the following results. There exists: 
\begin{enumerate}
\item an algebra  isomorphism $\varphi_h: \uqg \rightarrow U\g[[h]]$ and 
\item a `twist', \ie an element $\F\in U\g[[h]]\ot U\g[[h]]$ satisfying 
the relations
\eqa
&&(\varepsilon\ot \id)\F=\1=(\id\ot \varepsilon)\F
\label{cond2}\\
&&\F=\1\ot\1 +O(h)
\label{cond2bis}
\ena
($\1$ denotes the unit of $U\g$; (\ref{cond2bis}) 
implies that $\F$ is invertible as a formal power series in $h$)
\end{enumerate}
such that $H_h$ can be obtained from $H$ 
through the following equations. Let
$\F=\F^{(1)}\ot \F^{(2)}$, $\F^{-1}=\F^{-1(1)}\ot \F^{-1(2)}$,  
in a Sweedler's notation with {\it upper} indices; 
in the RHS a sum of many terms is implicitly understood, \eg 
$\sum_i \F^{(1)}_i\ot \F^{(2)}_i$. Then
\eqa
&& m_h=\varphi_h^{-1}\circ m \circ(\varphi_h\ot\varphi_h),
\label{defm}\\
&& \varepsilon_h:=\varepsilon\circ \varphi_h, \label{defc}\\
&&\Delta_h(x)=(\varphi_h^{-1}\ot\varphi_h^{-1})\big
[\F\Delta[\varphi_h(x)]\F^{-1}\big],\label{defd}\\
&&S_h(x)=\varphi_h^{-1}[\gamma^{-1} S [\varphi_h(x)]\gamma],
\label{def1}\\
&&\R=[\varphi_h^{-1}\ot\varphi_h^{-1}](\Fp q^{t\over 2}\F^{-1}).
\label{defR}
\ena
Here $t:=\Delta(\C)-\1\ot \C-\C\ot \1$ is 
the canonical invariant element 
in $U\g\ot U\g$ ($\C$
is the quadratic Casimir), the maps $m,\varepsilon,\Delta,S$ have
been linearly extended from $U\g$ to $U\g[[h]]$, and
\eq
\gamma := S \F^{-1(1)}\cdot\F^{-1(2)}, \z\x \qquad\qquad
\gamma^{-1} = \F^{(1)}\cdot S \F^{(2)}.
\label{def2}
\en

Equation (\ref{defm}) says that $U\g[[h]]$ and $\uqg$ are isomorphic 
through $\varphi_h$ as
algebras over $\cn[[h]]$. Equation (\ref{defc}) says that, up to this 
isomorphism, $\varepsilon_h$ and $\varepsilon$ coincide.
Equations (\ref{defd}), (\ref{def1}) say that, up to the same
isomorphism,
$\Delta_h,S_h$ differ from $\Delta,S$ by `similarity transformations'.
Equation (\ref{defd}) is not in contradiction with the coassociativity
of 
$\Delta$ and $\Delta_h$\footnote{To arrive at the above results
Drinfel'd
introduces the notion of quasitriangular quasi-Hopf algebra; 
the latter essentially involves the weakening of coassociativity
of the coproduct into a property  (``quasi-coassociativity'')
valid only up to a similarity 
transformation through an element $\phi\in U\g[[h]]^{\ot^3} $ 
(the ``coassociator''). This notion is useful
because quasitriangular quasi-Hopf algebra
are mapped into each other under twists (even if the latter 
is not trivial). 
As an intermediate result, he
shows that $U\g[[h]]$, beside the trivial quasitriangular 
quasi-Hopf structure 
$(U\g[[h]],m,\Delta,\varepsilon,S , 
\R\equiv \1^{\ot^2},\phi\equiv \1^{\ot^3})$,
has a non trivial one $(U\g[[h]],m,\Delta,\varepsilon,S , 
\R= q^{t\over 2},\phi\neq \1^{\ot^3})$.}, 
because the (nontrivial) coassociator
\eq
\phi:=[(\Delta\ot \id)(\F^{-1})](\F^{-1}\ot \1)(\1 \ot \F)
[(\id \ot \Delta)(\F)
\label{defphi}
\en
commutes with $\Delta^{(2)}(U\g)$ (we denote by
$\Delta^{(2)}$ the two-fold
coproduct),
\eq
[\phi,\Delta^{(2)} (U\g)]=0.
\label{ginv}
\en
>From the properties of $\phi$ it follows also that
$\gamma^{-1}\gamma'\in\mbox{Centre(U\g)}$, and 
$S \gamma=\gamma'^{-1}$.

The above formulae can be read also in the other direction
as giving a construction procedure of quasitriangular Hopf algebrae. 
They can be also applied to {\it triangular} deformations of $H$
[one needs only to set $t\equiv0$ in (\ref{defR})],
that are quasitriangular deformations with a trivial coassociator, 
$\phi=\1\ot\1\ot\1$. In fact the theorem cited above 
generalizes an older theorem \cite{dr1}, also by Drinfel'd.

 The twist $\F$ is 
defined (and unique) up to the transformation 
\eq
\F\rightarrow \F T,
\label{trasfo0}
\en 
where $T$ is a $\g$-invariant 
[\ie commuting with $\Delta(U\g)$] element
of $U\g[[h]]^{\ot^2}$ such that 
\eq
T=\1\ot\1+O(h),\z\z
(\varepsilon\ot\id) T=\1=(\id\ot\varepsilon) T.
\en 
Under this transformation
\eq
\phi\rightarrow
[(\Delta\ot \id)(T^{-1})](T^{-1}\ot \1)\phi(\1 \ot T)[(\id \ot 
\Delta)(T).
\label{trasfo}
\en
A function 
\eq
T=T\left(\1\ot {\cal C}_i, {\cal C}_i \ot \1, \Delta({\cal C}_i)\right)
\label{hil14}
\en
of the Casimirs ${\cal C}_i\in U\g$ 
of $U\g$ and of their coproducts clearly is \g-invariant. 
We find it plausible
that any \g-invariant $T$ must be of this form; although
we have found in the literature yet no proof of this conjecture,
in the sequel we assume that this is true.

We will often use a `tensor notation' for our formulae:
$\F_{12,3}=(\Delta\ot \id)\F_{12}$,
$\F_{123,4}=(\Delta^{(2)}\ot \id)\F_{12}$, and so on,
and definition (\ref{defphi}) reads
$\phi\equiv\phi_{123}=\F_{12,3}^{-1}\F_{12}^{-1}\F_{23}\F_{1,23}$.
$\phi$ satisfies the equations
\eq
\begin{array}{rcl}
q^{t_{13}+t_{23}\over 2} & = & \phi_{231}^{-1}q^{t_{13}\over 2} 
\phi_{132}q^{t_{23}\over 2} \phi_{123}^{-1}\cr
q^{t_{12}+t_{13}\over 2} & = & \phi_{312}q^{t_{13}\over 2} 
\phi_{213}^{-1}q^{t_{12}\over 2} \phi_{123}
\end{array}
\label{simpler}
\en
(they follow from $(\Delta_h\ot\id)\R=\R_{13}\R_{23}$,
$(\id \ot\Delta_h)\R=\R_{13}\R_{12}$).

While for the twist \F, apart from its existence, 
very little explicit knowledge is available,
Kohno \cite{ko} and Drinfel'd \cite{dr3} have proved that,
up to the transformation (\ref{trasfo}), $\phi$ is given by
\eq
\phi_m=\hat g^{-1}(x) \check g(x), \qquad\qquad\qquad 0<x<1,
\label{lilla}
\en  
where $\hat g, \check g(x)$ are $U\g[[h]]^{\ot^3}$-valued 
`analytic'\footnote{In the sense that the coefficients
$g_n(x)$ appearing in the expansion
$g(x)=\sum\limits_{n=0}^{\infty}g_n(x)h^n$ of $g$ in
$h$-powers are analytic functions of $x$ with values
in a finite-dimensional subspace  of $U\g^{\ot^3}$.} solutions 
of the first order linear differential equation
\eq
\frac{dg}{dx}=\hbar \left( \frac{t_{12}}x +\frac{t_{23}}
{x-1}\right)g, \qquad\qquad\qquad 0<x<1
\label{skzeq}
\en
($\hbar=\frac h{2\pi i}$)
with the following asymptotic behaviour near the poles:
\eq
\hat g \cdot x^{-\hbar t_{12}}\stackrel{x \rightarrow 0}
{\longrightarrow} \1^{\ot^3}
\qquad\qquad\qquad 
\check g \cdot (1-x)^{-\hbar t_{23}}\stackrel
{x \rightarrow 1}{\longrightarrow} \1^{\ot^3}.
\label{asymp}
\en
Using eq. (\ref{skzeq}) it is straightforward to verify 
that the RHS of eq. 
(\ref{lilla}) is indeed independent of $x$.
\footnote{Kohno and Drinfel'd proved that
$\phi$ can be obtained 
as the `monodromy' of a system of three first order linear
partial differential equations in three complex variables $z_i$
(the socalled universal Knizhnik-Zamolodchikov \cite{kz} equations),
with an $U\g^{\ot^3}$-valued unknown $f$, 
${\partial f\over \partial z_i}=\hbar\sum_{j\neq i}{t_{ij}\over
z_i-z_j}f$.
 The system can be reduced to the equation (\ref{skzeq}) 
exploiting its  invariance under linear 
tranformations $z_i\rightarrow az_i+b$. For a review of these
results see for instance Ref. \cite{chari}.} Using eq.
(\ref{asymp}) we can take the limit of eq. (\ref{lilla}):
\eq
\phi_m=\lim_{x_0\rightarrow 0^+}x_0^{-\hbar t_{12}} \check g(x_0).
\label{lilla2}
\en  
We can formally solve equations (\ref{skzeq}),
(\ref{asymp}) for $\check g$ by a path ordered integral,
\eq
\check g(x_0) = \lim_{y_0\rightarrow 0^+}\left\{
\vec{P}\exp\left[-\hbar\int\limits^{1-y_0}_{x_0} dx\left({t_{12}\over x}
+{t_{23}\over x-1}\right)\right] y_0^{\hbar t_{23}}\right\}
\en
($\vec{P}[A(x)B(y)]:=A(x)B(y)\vartheta(y-x)+B(y)A(x)\vartheta(x-y)$),
whence
\eq
\phi_m=\lim_{x_0,y_0\rightarrow 0^+}\left\{x_0^{-\hbar t_{12}}
\vec{P}\exp\left[-\hbar\int\limits^{1-y_0}_{x_0}dx\left({t_{12}\over x}
+{t_{23}\over x-1}\right)\right] y_0^{\hbar t_{23}}\right\}.
\label{integral}
\en
Note that $\phi_m=\1^{\ot^3}+ O(h^2)$.
We will say that the twist $\F$ is {\it `minimal'}  if the corresponding
$\phi$ (\ref{defphi}) is equal to $\phi_m$ or is trivial,
respectively in the case of $H_h=\uqg$ or $H_h$ is a triangular 
deformation of $U\g$.

The algebra isomorphism $\varphi_h:\uqg\rightarrow U\g[[h]]$
is defined up to an inner automorphism (a `similarity transformation')
of $U\g[[h]]$, 
\eq
\varphi_{h,v}(x):=v\varphi_h(x)v^{-1},
\label{simtra}
\en 
with $v=\1+O(h)\in U\g[[h]]$.  It is easy to check that 
Drinfel'd theorem \cite{dr3} remains true provided one replaces
$\F$ by $\F_v:=(v\ot v)\F\Delta(v^{-1})$ and all the 
objects derived from $\F$ correspondingly; in particular, it is
easy to check that the coassociator $\phi$ remains unchanged,
because it is $\g$-invariant
\eq
\phi_v=\Delta^{(2)} (v)\phi \Delta^{(2)} (v^{-1})=\phi.
\en
The freedom in choosing $\varphi_h$ (and $\F$) is 
usually eliminated or reduced if one requires it to satisfy additional 
properties, such as to lead to a specific $*$-structure 
for $\uqg$. 
The Lie algebra $g=sl(2)$ is the only $\g$ for which
explicit $\varphi_h$'s are known.

In the sequel we shall often use 
Sweedler's notations with {\it lower} indices for the
coproducts:
$\Delta(x)\equiv x_{(1)}\ot x_{(2)}$ for the 
cocommutative coproduct 
(in the RHS a sum 
$\sum_i x^i_{(1)}\ot x^i_{(2)}$
of many terms is implicitly understood),
$\Delta^{(n-1)}(x) \equiv x_{(1)}\ot\ldots\ot x_{(n)}$
for the $(n\!-\!1)$-fold cocommutative coproduct and 
$\Delta_h(x)\equiv x_{(\bar 1)}\ot x_{(\bar 2)}$
(with barred indices) for
the non-cocommutative one. 

\subsect{Deforming group-covariant Heisenberg algebrae}
\label{cha}

The generic undeformed Heisenberg algebra
is generated by the unit $\idA$ and
elements $a^+_i$ and $a^j$
satisfying the (anti)commutation relations
\eqa
&& a^i \,a^j      = \pm a^j\,a^i\nonumber\\
&& a^+_i\, a^+_j = \pm a^+_j\,a^+_i \label{ccr}\\
&& a^i \, a^+_j = \delta_j^i\idA \pm  a^+_j \,a^i 
\nonumber
\ena
(the $\pm$ sign refers to Weyl and Clifford algebras respectively).
$a^+_i,a^j$ transform under the action of $U\g$ according to some law
\begin{equation}
    x\trc a^+_i=\rho(x)^j_ia^+_j,\qquad\qquad
      x\trc a^i=\rho(Sx)_j^ia^j;                  \label{covl}
\end{equation}
here $x\in U\g$ and $\rho$ denotes
some matrix representation of \g. We shall
call the corresponding algebra $\A$. When $x\in\g$
the antipode reduces to $Sx=-x$.
Clearly $a^i$ belong to a representation 
of $U\g$ which is the contragradient $\rho^\vee = \rho^T \circ S $ 
(${}^T$ is the transpose) of the one of $a^+_i$, $\rho$. 
Because of the linearity of the transformation (\ref{covl}) we shall
also
say that $a^+_i,a^i$ are ``covariant'', or
``tensors'', under $\trc$. The action
$\trc$ is extended to products of the
generators using the standard rules of tensor product representations
(technically speaking, using the coproduct $\Delta$ of 
the universal enveloping algebra $U\g$, see formula (\ref{ext}) below), 
and then linearly to all of $\A$, $\trc:U\g\times\A\rightarrow\A$; 
this is possible because the action
of $\g$ is manifestly compatible with the commutation relations 
(\ref{ccr}),
and makes \A~ into a (left) module algebra of $(H,\trc)$.
In the sequel we shall denote by $\trc$ also its linear extension
to the corresponding algebrae of power series in $h$.

For suitable $\rho$ (specified below) \A~ admits a
deformation $\Aq$ with the same Poincar\'e series 
and the following features. $\Aq$ is generated by the unit
$\idAh$ and elements $\tilde A^+_i,\tilde A^i$
fulfilling deformed commutation relations
(DCR)  which can be put in  the form
\eqa
&& \tilde A^i\tilde A^j =\pm \P^{ji}_{hk}\tilde A^k 
\tilde A^h \label{gqcr0}\\
&& \tilde A^+_i\tilde A^+_j = \pm
\P_{ij}^{hk}\tilde A^+_h\tilde A^+_k 
\label{gqcr1}\\
&& \tilde A^i\tilde A^+_j =   \delta^i_j\idA \pm\,
\tilde P^F{}{ih}_{jk}\tilde A^+_h\tilde A^k 
\label{gqcr2}
\ena
(as before, the upper and lower sign
refer to Weyl and Clifford algebras respectively) and
transforms under the action $\trl$ of $U_h\g$ according to the law
\begin{equation}
x\,\trl\tilde A^+_i = \rho_h{}^j_i(x)A^+_j,
\qquad\qquad
x\,\trl\tilde A^i = \rho_h{}^i_j(S_h x)A^j.
\label{qtrans}
\end{equation}
Here $x\in U_h\g$, $\rho_h{}$ is the quantum group deformation of
$\rho$,
whereas $\P$ and $\tilde{P}^F$ are two suitable quantum-group-covariant
deformations of the ordinary permutator matrix $P$ 
(by definition $P^{ij}_{hk}=\delta^i_k\delta^j_h$). 
By a redefinition (\ref{simtra}) one can always choose a 
$\varphi_h$ such that
\eq
\rho_h=\rho\circ\varphi_h.
\en
$\tilde A^i$ belong to a representation of $U_h\g$ which is the 
quantum group contragradient $\rho_h^\vee = \rho_h^T \circ S_h $ of the 
one of $\tilde A^+_i$, $\rho_h$. 
Because of the linearity of the transformation (\ref{qtrans}) we shall
also
say that $\tilde A^+_i,\tilde A^i$ are ``covariant'', or
``tensors'', under $\trl$. The action $\trl$ is extended to products of
the
generators  by the formula
 \eq
x\,\trl(ab)=(x_{(\bar 1)}\trl a)(x_{(\bar 2)}\trl b)
\label{ext}
\en
and then linearly to all of $\Aq$, $\trl:U_h\g\times\Aq\rightarrow\Aq$; 
this is possible because the action
$\trl$ is compatible with the above commutation relations, 
and makes \Aq~ into a (left) module algebra of $(H_h,\trl)$.
The latter means also that 
\eq
(xy)\,\trl a=x\,\trl(y\,\trl a)
\label{modalg} 
\en
$\forall x,y\in\uqg$.
In the undeformed setting the formulae corresponding to the previous two 
are obtained by just replacing $\trl$ by $\trc$ and 
$\Delta_h(x)\equiv x_{(\bar 1)}\ot x_{(\bar 2)}$
by the cocommutative coproduct $\Delta(x)\equiv x_{(1)}\ot x_{(2)}$.
 
Up to our knowledge, only the following 
deformed algebras with the same Poincar\'e series as their
classical counterpart have been constructed:
${\cal A}^h_{+,sl(N),\rho_d}$,\cite{puwo,wezu}
${\cal A}^h_{-,sl(N),\rho_d}$,\cite{pusz}
${\cal A}^h_{+,so(N),\rho_d}$, \cite{cawa}, 
${\cal A}^h_{-,sp({N\over 2}),\rho_d}$ \cite{fionuovo}
(for each of the above Lie algebras
$\rho_d$ denotes the defining, $N$-dim representation)
as well as those $\Aq$ where $\rho$  
is the direct sum of many copies ($m$, say) of $\rho_d$,  
$\rho=\bigoplus\limits_{\mu=1}^m\rho_{d,\mu}$ \cite{fionuovo,quesne}
\footnote{For a generic $\rho$, second degree relations of the type
(\ref{gqcr0}-\ref{gqcr2}) may be not enough for a consistent definition
of \Aq; more precisely associativity may require the
introduction of third or higher degree relations, which have no
classical 
counterpart. In this case
the Poincar\'e series of \Aq~ will be smaller than
its classical counterpart, and $\Aq$ will be physically not so
interesting.}.

Explicitly, in the case $\rho=\rho_d$ the matrix $\tilde{P}^F$ is equal
to
$q\hat R$ (or to its inverse), where $\hat R$ is the braid matrix 
\cite{frt} of $U_h\g$,
\eq
\hat R:= c_{\small\g}\,P\,\left[(\rho_{d,h})^{\ot^2}(\R)\right] 
\qquad\qquad  c_{\small\g}:=
\cases{q^{\frac 1N}\x\mbox{if \g}=sl(N) \cr 1 \mbox{~~otherwise}\cr}
\label{defbraid}
\en
(here $\rho_{d,h}$ denotes the deformation of $\rho_d$
and we have introduced the factor $c_{\small\g}$ 
to match to the conventional normalization), whereas
$\P$ is a suitable first or second degree 
polynomial in $\hat R$. From formula (\ref{defR}) it then follows that
\eqa
\P & =& F\,U\,F^{-1} \label{genperm}\\
\tilde{\P} & =& F\,V\,F^{-1}
\label{genrmat}
\ena
where $F:=\rho_d^{\ot^2}(\F)$,
$V=c_{\small\g}\,P\,q^{\rho_d^{\ot^2}(t/2)}$
and $U$ is also a polynomial in $P\,q^{\rho_d^{\ot^2}(t/2)}$.
One may actually
choose $U=P$ without affecting the commutation relations, since in 
this case one can easily show that $(\1\mp U)\propto (1\mp P)$,
the (anti)symmetric projector.
For later use we
recall that from the projector decomposition and the
properties of $\hat R$  \cite{frt,wezu,cawa} it follows
that the `$q$-number operator' $\tilde {\cal N}:=\tilde A^+_i \tilde
A^i$
of ${\cal A}_{\pm,sl(N)\rho_d}^h$ satisfies the relations
\eq
\tilde {\cal N}\tilde A^+_i=\tilde A^+_i+q^{\pm 2} \tilde A^+_i 
\tilde{\cal N}
\z\z
\tilde{\cal N}\tilde A^i=q^{\mp 2}(-\tilde A^i+ \tilde A^i \tilde{\cal
N}),
\label{Ncr}
\en
and the invariant elements 
$\tilde A^+C\tilde A^+:=\tilde A^+_iC^{ij}\tilde A^+_j$,
$\tilde AC\tilde A:=\tilde A^iC_{ji}\tilde A^j$ 
of ${\cal A}_{+,so(N)\rho_d}^h$ satisfy the relations
\eqa
&&(\tilde A C\tilde A) \,\tilde A^i \:-\: \tilde A^i \, 
(\tilde A C\tilde A)   =   0 \label{pappa1} \\
&&(\tilde A^+ C\tilde A^+) \,\tilde A^+_i \:-\: \tilde A^+_i
(\tilde A^+ C\tilde A^+)   =   0 \label{pappa2} \\
&&(\tilde AC\tilde A)\,\tilde A^+_i \:-\: q^2\,
\tilde A^+_i\,(\tilde AC\tilde A)   =  
(1\!+\!q^{2-N})C_{ij}\tilde A^j\label{pappa3} \\
&&\tilde A^i\,(\tilde A^+C\tilde A^+)\:-\:
q^2\,(\tilde A^+C\tilde A^+)\,\tilde A^i 
  =   (1\!+\!q^{2-N})C^{ij}\tilde A^+_j.\label{pappa4}
\ena
In the case that $\rho$ is the direct sum of many copies (say $m$)
of $\rho_d$ the commutation relations
between the different copies are not
trivial. We refer the reader to Ref. \cite{fionuovo} for the explicit
form
of $\P,\tilde{\P}$. The important point here is that one can show 
that also in this case these matrices can be put in
the form (\ref{genperm}), (\ref{genrmat}), where 
\eq
F:=\rho^{\ot^2}(\F),
\label{Fmatrix}
\en
and $U,V$ are suitable $m\,N\,\times\, m\,N$ matrices such that
\eq
\left[U,\rho^{\ot^2}\left(\Delta(U\g)\right)\right]=
\left[V,\rho^{\ot^2}\left(\Delta(U\g)\right)\right]=0.
\label{UVmatrix}
\en

\renewcommand{\theequation}{\thesection.\arabic{equation}}
\sect{Realization of the quantum group action and of
candidates for the deformed generators}
\label{qcov}

It is immediate to check that one can define a Lie algebra homomorphism
$\sigma: \g\rightarrow \A$ by setting 
\eq 
 \sigma(x):=
\rho(x)^i_ja^+_i a^j
\label{jordan}
\en
for all $X\in \g$, and therefore extend it to
all of $U\g$ as an algebra
homomorphism $\sigma: U\g\rightarrow \A$ by setting
on the unit element $\sigma(\1_{U\g}):=\idA$. $\sigma$ can be seen
as the generalization of the Jordan-Schwinger
realization of $\g=su(2)$ \cite{bied}
\eq
\sigma(j_+)=a^+_{\up}a^{\downarrow},\qquad\qquad
\sigma(j_-)=a^+_{\downarrow}a^{\up},\qquad\qquad
\sigma(j_0)=\frac 12(a^+_{\up}a^{\up}-
a^+_{\downarrow}a^{\downarrow}).
\label{homo}
\en
Extending $\sigma$ linearly to the corresponding algebrae
of power series in $h$, we can define also an algebra
homomorphism 
\eq
\sigma_{\varphi_h}:=\sigma\circ\varphi_h:\uqg\rightarrow \A[[h]].
\en

Since we know that a deforming map exists, although we cannot write it 
explicitly we can say that it must be possible to construct the map 
$\tr$ defined
in the introduction. Our first step is to guess such a realization
of $\trl$ on $\Aq$. 
This requires fulfilling the conditions (\ref{ext}), (\ref{modalg}),  
which characterize a left module algebra.
There is a simple way to find such a realization, namely

\begin{prop} \cite{fio}
The (left) action $\tr:\uqg\times \A[[h]]\rightarrow \A[[h]]$ can be 
realized in an `adjoint-like' way:
\eq
x\tr a := \sigma_{\varphi_h}(x_{(\bar 1)}) a 
\sigma_{\varphi_h}(S_h x_{(\bar 2)}).
\label{defprop}
\en
\end{prop}
Using the basic axioms characterizing the coproduct, counit, 
antipode in a generic Hopf algebra
it is easy to check that (\ref{ext}), (\ref{modalg}) are indeed
fulfilled.
The realization (\ref{defprop}) 
is suggested by the cocommutative case, where it reduces to
\eq
x\trc a=\sigma(x_{(1)})\: a\: \sigma(S  x_{(2)}).
\label{cov}
\en 

Our second step is to realize elements $A^i,A^+_j\in\A[[h]]$ that
transform under (\ref{defprop}) as $\tilde A^i,\tilde A^+_j$
in (\ref{qtrans}). Note that $a^i,a^+_j$ do {\it not} transform
in this way. We recall 

\begin{prop} \cite{fio}
Let $\F$ be a twist associated to $\varphi_h$.
For any choice of $\g$-invariant elements $u,\hat u,v,\hat v=\idA+O(h)$ 
in $\A[[h]]$ (in particular if they are trivial) the elements
\eq
\begin{array}{lll}
A_i^+ &:= & \sigma(\F^{(1)})\, v\, a_i^+ \, \hat v\,
\sigma(S \F^{(2)}\gamma)\,  \cr
A^i&:= & \sigma(\gamma'S \F^{-1(2)})\, u\, a^i\,\hat u \, 
\sigma(\F^{-1(1)})         \cr                
\end{array}
\label{def3}
\en
transform under (\ref{defprop}) as $\tilde A^i,\tilde A^+_j$
in (\ref{qtrans}), and go
to $a^+_i,a^i$ in the limit $h\rightarrow 0$.\footnote{
The Ansatz (\ref{def3}) has some resemblance with
the one in Ref.  \cite{majid2}, prop. 3.3, which defines
an intertwiner $\alpha: U\g [[h]]\rightarrow \uqg$ of 
\uqg-modules. } 
\label{prop1}
\end{prop}
Looking for $\g$-invariants making $A^i,A^+_i$ fulfil the DCR
will be the third step of the construction (Sect. \ref{reqcr}).
 
{\it Remark 1.} 
Without loss of generality, we can assume 
$\F$ in definitions (\ref{def3}) to be minimal.
In fact, since any other twist can be written
in the form $\F T$ with $T\equiv T^{(1)}\ot T^{(2)}$ 
as in (\ref{hil14}), one finds that
\eqa
&&\sigma(T^{(1)}) (v\, a_i^+ \,\hat v)\sigma(S T^{(2)})=v'\, a_i^+ \,
\hat v' \\
&&\sigma(S T^{-1}{}^{(1)})(u\, a^i\,\hat u)
\sigma(T^{-1}{}^{(2)})=u'\, a^i\,\hat u' 
\label{change}
\ena
with $u',\hat u',v',\hat v'=\1+O(h)$ also $g$-invariant. This follows
from
the following observations: first, 
$\sigma({\cal C}_i),\sigma(S {\cal C}_i)$
are $\g$-invariant and central in $\sigma(U\g[[h]])$, therefore the
dependence of $T$ on ${\cal C}_i\ot \1$, $\1\ot{\cal C}_i$
translates just into
some replacements $u\rightarrow u'$, $v \rightarrow v'$, 
etc; second, $(v\, a_i^+ \,\hat v)$
transforms under $\trc$ exactly as $a_i^+$, whence
$\sigma({\cal C}_{i(1)})(v\, a_i^+ \,\hat v) \sigma(S {\cal C}_{i(2)})
\stackrel{(\ref{cov})}{=}
{\cal C}_i\trc(v\, a_i^+ \,\hat v) =c_i(v\, a_i^+ \,\hat v)$,
where $c_i\in\cn$ is the value
of the Casimir ${\cal C}_i$ in (the irreducible component of)
the representation $\rho$ to which $a^+_i$ belongs, and similarly for 
$(u\, a^i\,\hat u)$, in other words
the dependence of $T$ on ${\cal C}_{i(1)}\ot {\cal C}_{i(2)}$
translates just into mutiplication of the generators by a constant.

{\it Remark 2.} Note that if $\rho$ is reducible the previous
proposition
holds also if  we allow for {\it different} 
invariants $u,\hat u,v,\hat v$ within each irredicible component
of $\rho$.

In the sequel we shall often use the compact notation
\eq
\aa^+_i:=(v\, a_i^+ \,\hat v) \qquad \aa^i:=(u\, a^i\,\hat u).
\label{defaa}
\en

In the appendix we prove the following

\begin{lemma} 
If $\F$ is a `minimal', then
\eqa
\F &=& \gamma^{-1}(S \F^{-1(1)})
\F^{-1(2)}_{~(1)}\ot \F^{-1(2)}_{~(2)}\label{dritto1}\\
 & = & \F^{-1(1)}_{~(1)}\ot \gamma'(S \F^{-1(2)})
\F^{-1(2)}_{~(1)} \label{dritto2}\\
\F^{-1}&=& \F^{(1)}_{(1)}\ot \F^{(1)}_{(2)}
(S \F^{(2)})\gamma \label{dritto3}\\
 &= & \F^{(2)}_{(1)}
(S \F^{(1)})\gamma'^{-1}\ot \F^{(2)}_{(2)}. 
\label{dritto4}
\ena
\label{lemma2}
\end{lemma}
\vskip-1.0cm
We can find now useful alternative 
expressions for $A^+_i,A^i$. 
\begin{prop}
With a `minimal' $\F$, definitions (\ref{def3})  amount to
\eqa
A_i^+ &= &\aa^+_l\sigma(\F^{-1(2)})\rho(\F^{-1(1)})_i^l\label{def4.1}\\
A_i^+ &= &\rho(S \F^{(1)}\gamma'^{-1})_i^l\sigma(\F^{(2)})\aa^+_l
\label{def4.2}\\
A^i &= &\rho(\F^{(1)})^i_l\sigma(\F^{(2)}) \aa^l\label{def4.3} \\
A^i &= &\aa^l\sigma(\F^{-1(2)})\rho(\gamma^{-1} S \F^{-1(1)})^i_l.
\label{def4.4}
\ena
\label{prop2}
\end{prop}
\vskip-1.0cm
{\it Remark 3.} In spite of its original definition (\ref{defprop}),
from the latter expressions we realize that only a `semiuniversal form' 
of the type $(\rho\ot\id)\F^{\pm 1}$ for $\F$ is involved in the
definition of $A^i,A^+_j$.

$Proof~of~Prop.$ \ref{prop2}. Observing that
\eqa
\sigma(x)a&=&\sigma(x_{(1)} ) a\sigma(S x_{(2)}\cdot x_{(3)} )\\
a\sigma(x)&=&\sigma(x_{(3)} S  x_{(2)} )a\sigma(x_{(1)} )
\label{com}
\ena
for all $x\in U\g$, $a\in\A$, we find
\eqa
&& A_i^+ \stackrel{(\ref{cov}),(\ref{covl})}{=} 
\aa_l^+\sigma\left[\F^{(1)}_{(2)} 
(S \F^{(2)})\gamma\right]\rho(\F^{(1)}_{(1)})_i^l
\stackrel{(\ref{dritto3})}{=}
\aa_l^+\sigma(\F_1^{-1(2)})\rho(\F_1^{-1(1)})_i^l,\z\\
&& A^i \stackrel{(\ref{cov}),(\ref{covl})}{=} 
\rho(\F^{-1(1)}_{(1)} )^i_l\sigma\left[\gamma'(S \F^{-1(2)})
\F^{-1(1)}_{(2)}\right] \aa^l
\stackrel{(\ref{dritto2})}{=}\rho^i_l(\F^{(1)}_2)
\sigma(\F^{(2)}_2)\aa^l.\z
\label{def5}
\ena

Similarly one proves the other relations. 
$\Box$

\renewcommand{\theequation}{\thesection.\arabic{equation}}
\sect{Classical versus quantum group invariants}
\label{inva}

We have defined two actions $\trc,\tr$ on  $\A[[h]]$.
Their respective 
invariant subalgebras are respectively defined by
\eqa
\Ai[[h]] &: = &\{I\in \A[[h]]\x |\x x\trc I=\varepsilon(x) I\z\forall 
x\in U\g\}                      \label{def7} \\
\Aiq[[h]] &: = &\{I\in \A[[h]]\x |\x x\tr I=\varepsilon_h(x) I\z\forall 
x\in \uqg\}.                    \label{def7bis}
\ena
If $f$ is a deforming map corresponding to $\tr$,
the second subalgebra clearly contains $f(\Aiq)$, where $\Aiq$ denotes
the $U_h\g$-invariant subalgebra of $\Aq$.
What is the relation between $\Ai[[h]]$ and $\Aiq[[h]]$?

\begin{prop}
 $\Ai[[h]]=\Aiq[[h]]$.
\end{prop}
$Proof$. We show that both subalgebras coincide with the one
\eq
\{ I\in \A[[h]]\x |\x [\sigma(U\g), I]=0
\}.                    \label{nanna}
\en
Given any $I\in \Ai[[h]]$, $y\in U\g $, we find 
\eqa
I\sigma(y)&=&I\varepsilon(y_{(1)})\sigma( y_{(2)})
\stackrel{(\ref{def7})}{=}(y_{(1)}\trc I)\sigma( y_{(2)})\nonumber \\
&\stackrel{(\ref{cov})}{=}&\sigma(y_{(1)})I\sigma(S  y_{(2)}\cdot
y_{(3)})=\sigma(y_{(1)})I\varepsilon(y_{(2)})=\sigma(y)I.\nonumber
\ena
this proves that (\ref{nanna}) contains $\Ai[[h]]$. 
Conversely, if $[I,\sigma(U\g)]=0$ then
\[
x\trc I\stackrel{(\ref{cov})}{=}\sigma( x_{( 1)})I\sigma(S x_{(2)})
=I\sigma( x_{(1)}\cdot S x_{(2)})=I\varepsilon(x)
\]
for any $x\in U\g$,
proving that the set (\ref{nanna}) is contained in $\Ai[[h]]$. 
Replacing in the previous arguments 
$\sigma,\varepsilon,\trc$
and $\Delta(x)\equiv x_{(1)}\ot x_{(2)}$ by
$\sigma_{\varphi_h},\varepsilon_h,\tr$ and  
$\Delta_h(x)\equiv x_{(\bar 1)}\ot x_{(\bar 2)}$,
one proves that also $\Aiq[[h]]$ coincides with the algebra
(\ref{nanna}). 
$\Box$

In other words, the propostion
states that invariants under the $\g$-action $\trc$ are also
$U_h\g$-invariants (under $\tr$), and conversely, although in general
$\g$-covariant objects (\ie tensors) and $U_h\g$-covariant ones do not
coincide.

Since $\trc$ (resp. $\tr$) acts in a linear and homogeneous way
on the generators $a^i, a^+_j$ (resp. $A^i,A^+_j$), in the vector space 
$\Ai[[h]]=\Aiq[[h]]$ we can 
choose a basis $\{I^n \}_{n\in\nn}$ (resp. $\{I_h^n\}_{n\in\nn}$)
 consisting of normal
ordered homogeneous polynomials in 
$a^i, a^+_j$ (resp. $A^i,A^+_j$).  For all \g~ the simplest
invariant is the `number of particle operator' 
$n \equiv I^1:=a^+_ia^i$ and its deformed counterpart
${\cal N}\equiv I^1_h=A^+_iA^i$. In general the invariants will take the 
form
\eqa
I^n  &:=& a^+_{j_1}...a^+_{j_{k_n}}d^{j_1....j_{k_n}}_{i_1...i_{h_n}}
a^{i_1}...a^{i_{h_n}} \\
I_h^n &:=& A^+_{j_1}...A^+_{j_{k_n}}D^{j_1....j_{k_n}}_{i_1...i_{h_n}}
A^{i_1}...A^{i_{h_n}}
\label{def9}
\ena
($k_n,h_n\in\nn\cup\{0\}$); the coefficients 
$d^{j_1....j_{k_n}}_{i_1...i_{h_n}}$ 
(resp. $D^{j_1....j_{k_n}}_{i_1...i_{h_n}}$) 
depend on the particular $\g$ picked up and make up
classical (resp. quantum) group isotropic tensors, \ie satisfy
\eqa
\left[\left( \rho^{\ot^{k_n}}
\!\ot \!\rho^{\vee~\ot^{h_n}}\right)
\Big(\Delta^{(b_n\!-\! 1)}(x)\Big)
\right]^{J_nI_n}_{J_n'I'_n}d^{J'_n}_{I'_n}
&\!=\!& \varepsilon(x)d^{J_n}_{I_n}
\label{cisotrop} \\
\left[\left( \rho_h^{\ot^{k_n}}
\!\ot \!\rho_h^{\vee~\ot^{h_n}}\right)
\Big(\Delta_h^{(b_n\!-\! 1)}(y)\Big)
\right]^{J_nI_n}_{J_n'I'_n}D^{J'_n}_{I'_n}
&\!=\!& \varepsilon_h(y)D^{J_n}_{I_n}
\ena
$\forall x\in U\g[[h]], y\in\uqg$. Here and in the rest of the section
 we use the collective-index notation $ I_n\equiv(i_1....i_{h_n})$,
$ J_n\equiv(j_1....j_{k_n})$ and   the
short-hand notation $b_n\!:=\!h_n\!+\!k_n$.
Using formula (\ref{def1}) it is 
straightforward  to verify that the $d$'s 
and $D$'s are related to each other by
\eq
 D^{ J_n}_{ I_n}\propto
\left[\left( \rho^{\ot^{k_n}}\!\ot \!\rho^{\vee~\ot^{h_n}}\right)
\Big((\1^{\ot^{k_n}}\ot (\gamma'^{-1})^{\ot^{h_n}})
\F_{12....b_n}\Big)\right]^{ J_n I_n}_{ J'_n I'_n}
d^{ J'_n}_{ I'_n},
\label{propto}
\en
where $\F_{12....b}\in U\g[[h]]^{\ot^b}$  
is an intertwiner between $\Delta^{(b-1)}$ and
$\Delta_h^{(b-1)}$ and is given,
up to multiplication from the right by a \g-invariant tensor
$Q\in U\g[[h]]^{\ot^{b_n}}$, by
\eq
\F_{12....b}=\F_{(b\!-\!1)b}\F_{b\!-\!2,(b\!-\!1)b}....
\F_{1,2...b}.
\label{def10}
\en
The isotropic tensors corresponding to $I^1,I^1_h$ are
$d^j_i=\delta^j_i=D^i_j$.
The replacement $\F\rightarrow\F\cdot T$, with 
$T\in U\g[[h]]^{\ot^2}$ and \g-invariant, results also
in multiplication from the right by a related $Q$.

Relation (\ref{propto}) guarantees the existence of $D$'s in
one-to-one correspondence with the $d$'s, but from the 
practical viewpoint is not of much help for finding the $D$'s
(since the universal $\F$ is unknown and its matrix
representations are known only for few representations);
the latter can be found more easily from the knowledge
of \R and a direct study of $\trl$.

The question whether 
$\Ai$ is finitely generated, \ie whether all $I^n$ can be
expressed as {\it polynomials} in a finite number of them,
is part of an important problem originally raised by 
Hilbert\footnote{The fourteenth problem proposed by
Hilbert at the 1900 conference of mathematicians.}.
This is in general not the case. However, they can
be expressed as {\it algebraic functions} of 
a finite number of them \cite{dieu}. The classification
of the latter (for arbitrary \g)
is not completed yet ~\cite{azc}. 

We would like to ask here a different question. 
According to the above proposition, $I^n_h$
can be expressed as a power series in $h$ with coefficients in
$\Ai$, in other words as a `function' $I^n_h=m^n(h,\{I^n \})$.
Since the latter are functions of $a^i,a^+_j$,
we can express $I^n_h$ also as a function 
$I^n_h=f^n(h,a^i,a^+_j)$. How to find $m^n, f^n$?
We give an explicit
answer to the second question by proving (see the appendix)

\begin{prop}
\eqa
&& I_h^n  =  (\aa^+...\aa^+)_{ M_n}(\aa...\aa)^{ L_n}
\times \label{formulinv}\cr
& &\left[\left( \rho^{\ot^{k_n}}
\!\ot \!\rho^{\vee~\ot^{h_n}}\ot\sigma\right)
\left(\phi^{-1}_{(b_n\!-\!1)b_n(b_n\!+\!1)}
\phi^{-1}_{(b_n\!-\!2),(b_n\!-\!1)b_n,(b_n\!+\!1)}
...\phi^{-1}_{1,2...b_n,(b_n\!+\!1)}\right)
\right]^{ M_n L^n}_{ J_n I_n}
d^{ J_n}_{ I_n}
\nonumber
\ena
where 
$\phi_{1,2...m,m\!+\!1}:=(\id\ot\Delta^{(m-2)}\ot\id)\phi_{123}$
and $b_n\!:=\!h_n\!+\!k_n$.
\label{propinv}
\end{prop}

{\it Remark 4.} Note that  in these equations 
the whole dependence on the twist $\F$ is concentrated in
the coassociator $\phi$ of \g and in its coproducts. Consequently,
use of formula (\ref{integral}) allows the explicit determination of
the dependence of $I_h^n$'s on $a^i, a^+_j$; clearly,
the latter will be in general highly non-polynomial.
 
Let us address now the first question.
If $H_h$ is {\it triangular} then $\phi^{-1}$ and all
its coproducts are trivial, the $\g$-invariants
$u,v,\hat u,\hat v$ may be chosen trivial as well \cite{fio},
and from the previous proposition we find
$I_h^n=I^n$.
If $H_h$ is a {\it genuine quasitriangular} Hopf algebra, such as
$\uqg$, then in general
$I_h^n\neq I^n$;
the $I_h^n$ will be some nontrivial function of the 
$I^m $'s, generally speaking highly non-polynomial, as well.
This can be already verified for the simplest invariants.
We will show in next section that \eg $I^1_h=(n)_{q^2}\equiv
(I^1)_{q^2}$
in the $\g=sl(N)$ case. In general,
\eq
I^n_h=I^n+O(h).
\en
This follows from $D^{I_n}_{J_n}= d^{I_n}_{J_n}+O(h)$,
$A^i=a^i+O(h)$, $A^+_i=a^+_i+O(h)$. 

In the $\g=so(N),\rho=\rho_d$ case, beside $\delta^i_j$,
another basic isotropic tensor is the classical metric
matrix $c_{ij}=c_{ji}$ (with inverse $c^{ij}=c^{ji}$, 
to which there corresponds the deformed
metric matrix \cite{frt} $C_{ij}$, and its inverse
$C^{ij}$): 
\eq
C^{ij}=F^{ij}_{hk}c^{hk}=\rho_d(\gamma^{-1})^i_hc^{hj}.
\z\z C_{ij}=c_{ih}\rho_d(\gamma)^h_j;
\label{metmat}
\en
the last two equalities follow from the $so(N)$ property
\eq
\rho_d(S x)^i_j=\rho_d(x)^m_lc^{li}c_{mj} \z\z x\in U\g.
\label{partico}
\en
So one can build the invariants
\eq
\begin{array}{rclcrcl}
I^{2,0}  &:=& a^ic_{ij}a^j\equiv a\,c\,a &\z &
I^{0,2}  &:=& a^+_ic^{ij}a^+_j\equiv a^+\,c\,a^+ \cr
I_h^{2,0} &:=& A^iC_{ji}A^j\equiv A\,C\,A &\z &
I_h^{0,2} &:=& A^+_iC^{ij}A^+_j\equiv A^+\,C\,A^+; \cr
\end{array}
\label{scal}
\en
we will see in next section that $I_h^{2,0}\neq I^{2,0} $,
$I_h^{0,2}\neq I^{0,2} $.
We leave the determination of the general dependence of $I^n_h$ on
$I^m$'s as a subject for further investigation.

\renewcommand{\theequation}{\thesection.\arabic{equation}}
\renewcommand{\theequation}{\thesubsection.\arabic{equation}}
\sect{Realization of the deformed generators}
\label{reqcr}

In section \ref{qcov} we have left some freedom
in the definition of $A^i,A^+_i$:
the \g-invariants $u,\hat u,v,\hat v$ appearing in the definitions
(\ref{defaa}) of $\aa^i,\aa^+_i$ have not been specified.
Can we choose them in such a way that $A^i,A^+_i$
fulfil the DCR (deformed commutation relations) of $\Aq$? To 
explicitly study
this question in the appendix we prove 

\begin{prop}
If we replace $\tilde A^i,\tilde A^+_j\rightarrow A^i,A^+_j$
[with $A^i,A^+_j$ defined as in
formulae (\ref{def4.1}), (\ref{def4.3}), with a minimal $\F$], then
equations (\ref{gqcr0}), (\ref{gqcr1}), (\ref{gqcr2})
become equivalent to
\eqa
\aa^i\, \aa^j & = & \:\pm\: (M^{-1}\,U\, M)_{lm}^{ji}\aa^m\, \aa^l
\label{figata1}\\
\aa^+_i\,\aa^+_j\,& = & \:\pm\:
\aa^+_l\,\aa^+_m\,(M^{-1}\,U\,M)^{lm}_{ij} 
\label{figata2}\\
 \aa^i\,\aa^+_j & = & \delta^i_j\idA\:\pm\: 
\aa^+_l(M^{-1}\,V\,M)^{il}_{jm}\aa^m
\label{figata3}
\ena
where $U\equiv\Vert U^{ij}_{hk}\Vert$, $V\equiv\Vert V^{ij}_{hk}\Vert$ 
are the (numerical) matrices introduced in equations (\ref{UVmatrix})
 and $M\equiv\Vert M^{ij}_{hk}\Vert$
is the $\sigma(U\g[[h]])$-valued matrix defined by
\eq
M:= (\rho\ot \rho\ot \sigma)(\phi_m).
\label{def8}
\en
\label{figate}
\end{prop}
\vskip-1.0cm
 We recall that, if $\rho=\rho_d$, $U$ is the permutation
matrix $P$ and 
$V\propto P\,q^{\rho_d^{\ot^2} (\frac t2)}$.

{\it Remark 5.} The above equations have to be understood as equations
in the unknowns $u,\hat u,v,\hat v$. They can be studied explicitly
because the whole dependence  on $\F$  is concentrated again in
the coassociator $\phi$ of \uqg.

{\it Remark 6.}
If $H_h$ is a triangular deformation, then $U=V=P$,  $\phi=\1^{\ot^3}$ 
(and consequently $M=\1^{\ot^3}$), and the eq.
(\ref{figata1}), (\ref{figata2}) are
satisfied with trivial invariants $u,\hat u,v,\hat v$, \ie
with $\aa^i=a^i$, $\aa^+_i=a^+_i$. This was already shown in
Ref. \cite{fio}.

To look for solutions of eq. (\ref{figata1}), (\ref{figata2}), 
(\ref{figata3}) for genuine quasitriangular
deformations we have to treat the \g's belonging to different
classical series separatly.
We consider here ${\cal A}^h_{\pm,sl(N),\rho_d}$,
${\cal A}^h_{+,so(N),\rho_d}$.

\subsect{The case of ${\cal A}^h_{\pm,sl(N),\rho_d}$}
\label{better}

As a basis of \g we choose $\{E_{ij}\}_{i,j=1,...,N}$ with
$\sum\limits_{i=1}^NE_{ii}\equiv 0$ (so that there exist
only $N^2-1$ linearly independent $E_{ij}$), satisfying
\eq
[E_{ij},E_{hk}]=E_{ik}\delta_{jh}-E_{jh}\delta_{ik}
\en
The quadratic Casimir reads
\eq
{\cal C}= E_{ij}E_{ji},
\en
implying
\eq
t=2\, E_{ij}\ot E_{ji}
\en
The matrix representation of $E_{ij}$ in the fundamental
representation $\rho$ takes the form
\eq
\rho(E_{ij})=e_{ij}-\frac{\delta_{ij}}N \1_N,
\en
where $e_{ij}$ is the $N\times N$ matrix with
all vanishing entries but a 1 in the $i$-th row and $j$-th
column, and $\1_N=\sum_i e_{ii}$ is the $N\times N$ 
unit matrix; whereas
the Jordan-Schwinger realization takes the form
\eq
\sigma(E_{ij}) = a^+_ia^j-\frac{\delta_{ij}}N n.
\label{dindon}
\en
As a consequence
$\sigma({\cal C})=n(N\pm n\mp 1)-\frac{n^2}N$.
>From the previous equations one finds
\eqa
&&(\rho\ot\rho\ot\sigma) \left(\frac{t_{12}}2\right) =\,e_{ij}\ot e_{ji}
\ot\idA-\frac 1N
\1_N\ot\1_N\ot\idA\,=:\, P-\frac {\1^{\ot^3}}N, \nonumber\cr
&&(\rho\ot\rho\ot\sigma)\left(\frac{t_{23}}2\right) = \1_N\ot 
e_{ij}\ot a^+_ja^i-
\1_N\ot\1_N\ot \frac nN\,=:\,A-\1_N^{\ot^2}\ot \frac nN\label{def14}\cr
&&(\rho\ot\rho\ot\sigma)\left(\frac{t_{13}}2\right) = e_{ij}\ot 
\1_N\ot a^+_ja^i-
\1_N\ot\1_N\ot \frac nN\,=:\,B-\1_N^{\ot^2}\ot \frac nN ;\nonumber
\ena
$P$ denotes the permutation matrix on $\cn^N\ot \cn^N$, 
multiplied by $\idA$.

$n:=a^+_ia^i$ is an element of ${\cal A}^{inv}_{\pm,sl(N),\rho_d}[[h]]$. 
We try to solve eq. (\ref{figata1}-\ref{figata3})
with invariants $u,\hat u,v,\hat v$
depending only on $n$. Using relations 
\eq
[n,a^+_i]\,=\,a^+_i \z\z [n,a^i]\,=\,-\, a^i.
\en 
we can thus commute $\hat v$
to the left of $a^i$ and $u$ to the right of 
$a^+_i$ in formula (\ref{defaa}), and look for 
$\aa^i,\aa^+_i$ directly in the form
$\aa^i:=I a^i$, $\aa^+_i:=a^+_i\tilde I$, with
$I=I(n)\in\Ai[[h]]$, $\tilde I=\tilde I(n)\in\Ai[[h]]$.
>From eq. (\ref{def4.1}),  (\ref{def4.3})
it follows ${\cal N}:=A^+_iA^i=\aa^+_i\aa^i=n\,\hat I(n\!-\!1)$, where
$\hat I(n):=I(n)\tilde I(n)$. In order that ${\cal N},A^+_i,A^i$ 
satisfies the
commutations relations (\ref{Ncr}), we therefore require 
$\hat I(n)=\frac{(n\!+\!1)_{q^{\pm 2}}}{n\!+\!1}$,
with $(x)_a:={a^x-1\over a-1}$. 
Summing up, we leave $I(n)$ undetermined and we pick
\eqa
\aa^i &:=& I a^i,\z\z I\equiv I(n),\label{def13} \\
\aa^+_i &:=& a^+_i\tilde I\z\z \tilde I(n):=I^{-1}(n)\,
{(n\!+\!1)_{q^{\pm 2}}\over 
(n\!+\!1)}. \nonumber
\ena
These ansatz can also be written in the equivalent form
\eq
\aa^i=u(n)a^iu^{-1}(n), \z\z
\aa^+_i=v(n) a^+_iv^{-1}(n),
\label{ansatz0}
\en
where $u,v$ are constrained by the relation
\eq
u\,v^{-1}=y=y_{sl(N)}:=\frac{\Gamma(n+1)}
{\Gamma_{q^2}(n+1)}
\label{yslN}
\en
and $\Gamma,\Gamma_{q^2}$ are the Euler's $\Gamma$-function
and its deformation (\ref{defgammaq}).

We have now the right ansatz to show that the DCR of 
$N$-dimensional $U_hsl(N)$-covariant Heisenberg
algebra are fulfilled. In the appendix we prove

\begin{theorem}
When $\g=sl(N)$, the objects $A^i,A^+_i$ ($i=1,2,...,N$) 
defined in formulae
(\ref{def4.1}), (\ref{def4.3}), (\ref{def13}) 
satisfy the
corresponding DCR (\ref{gqcr0}), (\ref{gqcr1}), (\ref{gqcr2}).
\label{propsln}
\end{theorem}
\vskip-0.5cm
In Ref. \cite{fio} the case $\g=sl(2)$ was worked out explicitly.
Choosing $u=v^{-1}=\sqrt{y_{sl(2)}}$, we found for the 
$A^i,A^+_i\in{\cal A}_{+,sl(2),\rho_d}[[h]]$ ($i=\uparrow,\downarrow$),
\eq
\begin{array}{rclcrcl}
A^+_{\up} & = &\sqrt{(n^{\up})_{q^2}\over n^{\up}}
q^{n^{\downarrow}}a^+_{\up} &\qquad
\qquad A^+_{\downarrow} & = &
\sqrt{(n^{\downarrow})_{q^2}\over n^{\downarrow}}
a^+_{\downarrow}  \nonumber \\
A^{\up} & = & a^{\up}\sqrt{(n^{\up})_{q^2}\over n^{\up}}
q^{n^{\downarrow}} &\qquad
\qquad A^{\downarrow} & = &a^{\downarrow} 
\sqrt{(n^{\downarrow})_{q^2}\over n^{\downarrow}},
\end{array}
\label{lastb}
\en
and for the $A^i,A^+_i\in{\cal A}_{-,sl(2),\rho_d}[[h]]$
\eq
\begin{array}{rclcrcl}
A^+_{\up} & = & q^{-n^{\downarrow}}a^+_{\up} &\qquad
\qquad A^+_{\downarrow} & = &
a^+_{\downarrow}  \nonumber \\
A^{\up} & = &a^{\up}q^{-n^{\downarrow}} &\qquad
\qquad A^{\downarrow} & = & a^{\downarrow}.
\end{array}
\label{lastf}
\en
Here $n_i:=a^+_ia^i$ (no sum over $i$).
Let us compare the generators (\ref{lastb})  with the
ones found in Ref. \cite{oleg}. In our
notation the latter would read
\eq
\begin{array}{rclcrcl}
A^+_{\alpha\up} & = & q^{n ^{\downarrow}}a^+_{\up} &\qquad
\qquad A^+_{\alpha\downarrow} & = &
a^+_{\downarrow}  \nonumber \\
A^{\alpha\up} & = &a ^{\up}{(n ^{\up})_{q^2}\over n ^{\up}}
q^{n ^{\downarrow}} &\qquad
\qquad A^{\alpha\downarrow} & = &a ^{\downarrow} 
{(n ^{\downarrow})_{q^2}\over n ^{\downarrow}}.
\end{array}
\en
It is straightforward to check that
the element $\alpha\in\A[[h]]$
such that $A^{\alpha i}=\alpha A^i\alpha^{-1}$,
 $A^+_{\alpha i}=\alpha A^+_i\alpha^{-1}$
[formulae (\ref{innerauto})] is
\eq
\alpha:=\sqrt{\frac{\Gamma(n^{\up}+1)\Gamma(n^{\downarrow}+1)}
{\Gamma_{q^2}(n^{\up}+1)\Gamma_{q^2}(n^{\downarrow}+1)}},
\en
where $\Gamma$ is the Euler $\Gamma$-function and 
$\Gamma_{q^2}$ its $q$-deformation (\ref{defgammaq}).

\subsect{The case of ${\cal A}^h_{+,so(N),\rho_d}$}

As a basis of \g=$so(N)$ we choose $\{L_{ij}\}_{i,j=1,...,N}$ with
$L_{ij}=-L_{ji}$ (so that there exist
only $\frac {N(N\!-\!1)}2$ linearly independent $L_{ij}$), 
satisfying
\eq
[L_{ij},L_{hk}]=L_{ik}c_{jh}+L_{kj}c_{ih} -L_{hj}c_{ik}
-L_{ih}c_{jk};
\en
here $c_{ij}$ denotes the (classical) metric 
matrix on the
$N$-dimensional Euclidean space ($c_{ij}=c_{ji}$),
 which in the special
case we choose real Cartesian
coordinates takes simply the form $c_{ij}=\delta_{ij}$.
In the rest of this subsection classically-covariant 
indices will be lowered and raised by means of multiplication
by $c$: $v_i=c_{ij}v^j$ $v^i=c^{ij}v_j$, etc., and 
$v\cdot w:=v^iw^jc_{ij}=v_iw^i=v^iw_i$.
The quadratic Casimir reads
\eq
{\cal C}=\frac 12 L_{ij}L^{ji},
\en
implying
\eq
t= L_{ij}\ot L^{ji}
\en
The matrix representation of $E_{ij}$ in the fundamental
representation $\rho$ takes the form
\eq
\rho(L_{ij})=e_{ih}c_{hj}-e_{jh}c_{hi},
\en
and the Jordan-Schwinger realization becomes
\eq
l_{ij}:=\sigma(L_{ij})= a^+_ia^hc_{hj}-a^+_ja^hc_{hi}=
a^ha^+_ic_{hj}-a^ha^+_jc_{hi}.
\en

It is easy to work out 
\eq
l^2:=\sigma({\cal C})+(1\!-\!\frac N2)^2=\left(n+\frac 
N2-1\right)^2- (a^+\!\cdot\! a^+)(a\!\cdot\! a),
\label{def11}
\en
and to check that, as expected
\eq
[l^2,a^+\!\cdot\! a^+]=0=[l^2, a\cdot a]
\en
A direct calculation also shows that
\eqa
&&[l^2,a^i]  = -a^i(2n+1+N)+2(\!a\cdot a\!)a^+_jc^{ji} = 
-a^i(2n\!-\!3\!+\!N)+2a^+_jc^{ji}(\!a\cdot a\!)    \nonumber \cr
&&[l^2,a^+_i] =  a^+_i(2n+3+N)-2c_{ij}a^j(\!a^+\!\cdot\! a^+\!)
 =  a^+_i(2n\!-\!1\!+\!N)-2c_{ij}(\!a^+\!\cdot\! a^+\!)a^j.\nonumber
\ena
We look for ``eigenvectors'' of $l^2$ 
\[
l^2\alpha^i=\alpha^i\lambda \z\z
l^2\alpha^+_i=\alpha^+_i\mu, 
\]
in the form 
$\alpha^i=a^i\gamma+ a^+_jc^{ji}\, (a\!\cdot\! a)\,\delta$, 
$\alpha^+_i=a^+_i\alpha+ a^jc_{ji}\, (a^+\!\cdot\! a^+)\,\beta$
with ``eigenvalues'' $\lambda,\mu$ and ``coefficients''
$\alpha,\beta,\gamma,\delta$ depending on $n,l^2$. 
We find second order equations for $\lambda,\mu$ 
with solutions $\lambda,\mu=(l\pm 1)^2$, where
formally $l=\sqrt{l^2}$.  We can therefore
consistently extend ${\cal A}_{+,so(N),\rho_d}$ 
by the introduction of
a new generator $l$ [whose square is constrained to
give the $l^2$ defined in eq. (\ref{def11})]  such that
\eqa
&&\alpha^i_{\pm} := a^i(n\!+\!\frac N2 \!-\!1\!\pm l)
- c^{ij}a^+_j\,(a\cdot a) =  a^i(n\!+\!\frac N2 
\!+\!1\!\pm l) - (a\cdot a) \,c^{ij}a^+_j,  \nonumber\cr
&&\alpha^+_{i,\pm} :=  a^+_i(n\!+\!\frac N2\! -\!1\!\pm l)
- (a^+\!\cdot\! a^+)\,c_{ij}a^j  =  a^+_i(n\!+\!\frac N2 
\!+\!1\!\pm l) - c_{ij}a^j \,(a^+\!\cdot\! a^+)\nonumber
\ena
satisfy
\eqa
l\,\alpha^+_{i,\pm} & = & \alpha^+_{i,\pm}\,(l\!\pm\! 1),
\nonumber\\
l\,\alpha^i_{\pm} & = & \alpha^i_{\pm}\,(l\!\mp\! 1).
\label{eige}
\ena

After these preliminaries, let us determine
the right $\aa^i,\aa^+_i$'s for $A^i,A^+_i$
to satisfy the DCR. To satisfy at once
eq.'s (\ref{pappa1}),(\ref{pappa2})  
we make the ansatz: 
\eq
\aa^i=u(n,l)a^iu^{-1}(n,l), \z\z
\aa^+_i=v(n,l)a^+_iv^{-1}(n,l)
\label{ansatz}
\en
This implies 
\eqa
A^+\,C\, A^+&\stackrel{(\ref{ansatz})}{=}&
\aa^+_l\aa^+_mc^{lm}=v\, a^+\cdot a^+\, v^{-1}\\
A\,C\, A&\stackrel{(\ref{ansatz}}{=}&
\aa^l\aa^m c_{lm}=u\, a\cdot a\, u^{-1}
\ena

The DCR determine only the product $y:=v^{-1}u$;
we are going to show now that eq.'s
(\ref{pappa3}), (\ref{pappa4}) completely 
determine the latter.
It is immediate to check that the former implies
\[
 y\left[2c^{ij}a^+_j
(a^+\cdot a^+)a^i\right]y^{-1}(n\!+\!2,l)-q^2
y(n\!-\!2,l)(a^+\cdot a^+)a^iy^{-1}=(1\!+\!q^{N-2})
c^{ij}a^+_j.
\]
Expressing $a^+_i, c_{ij}(a^+\cdot a^+)a^j$ 
as  combinations of
$\alpha_{i\pm}$ we easily move $y$  past
the ``eigenvectors'' $\alpha_{i\pm}$ of $n,l$;
factoring out  (from the right)
${c_{ij}\over 2l}$ we end up with a LHS being a
 combination of 
$\alpha_{i,+}$, $\alpha_{i,-}$. Therefore eq.
(\ref{pappa3}) amounts to the condition that the
corresponding coefficients vanish:
\[
(1\!+\!q^{N-2})=(n\!+\!\frac{N}2\!+\!1\!-\!l)\,
y(n\!+\!1,l\!+\!1)y^{-1}(n\!+\!2,l)
-q^2(n\!+\!\frac{N}2\!-\!1\!-\!l)\,
y(n\!-\!1,l\!-\!1)y^{-1}(n,l)
\]
\[
(1\!+\!q^{N-2})=(n\!+\!\frac{N}2\!+\!1\!+\!l)\,
y(n\!+\!1,l\!+\!1)y^{-1}(n\!+\!2,l)
-q^2(n\!+\!\frac{N}2\!-\!1\!-\!l)\,
y(n\!-\!1,l\!-\!1)y^{-1}(n,l) 
\]
Similarly, from eq. (\ref{pappa4}) it follows
\[
(1\!+\!q^{N-2}) = (n\!+\!\frac{N}2\!+\!1\!-\!l)\,
y(n,l)y^{-1}(n\!+\!1,l\!-\!1)
-q^2(n\!+\!\frac{N}2\!-\!1\!-\!l)\,
y(n\!-\!2,l)y^{-1}(n\!-\!1,l\!-\!1)
\]
\[
(1\!+\!q^{N-2}) = (n\!+\!\frac{N}2\!+\!1\!+\!l)\,
y(n,l)y^{-1}(n\!+\!1,l\!+\!1)
-q^2(n\!+\!\frac{N}2\!-\!1\!+\!l)\,
y(n\!-\!2,l)y^{-1}(n\!-\!1,l\!+\!1)
\]
It is straightforward to check that the last four
equations  are solved by
\eq
u\,v^{-1}=
y=y_{so(N)}:=\left(\frac{1\!+\!q^{N-2}}2\right)^{-n}
\frac{\Gamma\left[\frac 12\left(n\!+\!\frac{N}2\!+\!1\!-\!l
\right)\right]\,\Gamma\left[\frac 12\left(n\!+\!
\frac{N}2\!+\!1\!+\!l\right)\right]}{\Gamma_{q^2}
\left[\frac 12\left(n\!+\!\frac{N}2\!+\!1\!-\!l
\right)\right]\,\Gamma_{q^2}\left[\frac 12
\left(n\!+\!\frac{N}2\!+\!1\!+\!l\right)\right]},
\label{solve}
\en
where $\Gamma,\Gamma_{q^2}$ are the Euler's $\Gamma$-function
and its deformation (\ref{defgammaq}).

We have now the right ansatz to fulfill  the DCR of 
$N$-dimensional $U_hso(N)$-covariant Weyl
algebra. We state without proof the

\begin{theorem}
When $\g=so(N)$, the objects $A^i,A^+_i$ ($i=1,2,...,N$) 
defined in formulae
(\ref{def4.1}), (\ref{def4.3}), (\ref{ansatz}), (\ref{solve})
satisfy the
corresponding DCR (\ref{gqcr0}), (\ref{gqcr1}), (\ref{gqcr2}).
\label{propson}
\end{theorem}

\renewcommand{\theequation}{\thesection.\arabic{equation}}
\sect{$*$-Structures}
\label{star}

Given the Hopf $*$-algebra
$H_h=(\uqg, m,\Delta_h,\varepsilon_h,S_h,\R,*_h)$,
we ask now whether the $*$-structures 
$\dagger_h$ of \Aq compatible with
the action $\tr$ of \uqg,\ie such that 
\eq
(x\,\tilde\tr a)^{\dagger_h}=S_h^{-1}(x^{*_h})\tilde\tr a^{\dagger_h},
\label{condstar}
\en
can be naturally realized by the ones of  \A.

We stick to the case that $H_h$ is the compact real section
of \uqg. Then \uqg as an algebra 
is isomorphic to $U\hat{\g}[[h]]$,
where $\hat{ \g}$ is the compact section of $\g$ and
$h\in\rn$, and the trivializing maps
$\varphi_h$ intertwine between ${*_h}$ and $*$,
$[\varphi_h(x)]^{* }=\varphi_h(x^{*_h})$
where $* $ is the classical $*$-structure
in $U\g$ having the elements of $\hat {g}$ as
fixed points.

If $\F$ is unitary then the corresponding 
$\gamma,\gamma',\phi$ clearly satisfy
\eq
\gamma'=\gamma^{* }\z\z\phi^{~{* }\ot {* }\ot {* }}
=\phi^{-1}.
\en
On the other hand, it is evident that the `minimal'
coassociator $\phi_m$ (\ref{lilla}) is also unitary
(because  $h\in\rn$); one could actually show that
the unitary $\F$ is also minimal.

If $\rho_h$ is a $*$-representation of
H, the $*$-structure
$(\tilde A^i)^{\dagger_h}=\tilde  A^+_i$
is clearly compatible with $\trl$
[condition (\ref{condstar})]; 
the classical counterpart of $\rho_h$ is also 
a $* $-representation
$\rho$ of $H$
(\ie $\rho(x^{* })=\overline{\rho^T(x)}$),
and formula
\eq
(a^i)^{\dagger }= a^+_i
\label{star1}
\en
defines in \A~ 
a $*$-structure (`hermitean conjugation')
${}^{\dagger }$ compatible with $\trc$. 
Correspondingly,
it is immediate to check
that $\sigma$, $\sigma_{\varphi_h}$
become $* ,{*_h}$-homomorphisms
respectively, 
\eq
\sigma\circ * =\dagger \circ \sigma \z\z
\sigma_{\varphi_h}\circ {*_h}=\dagger \circ \sigma_{\varphi_h},
\label{starhom}
\en
and $\tr$ as defined in formula  (\ref{defprop})  also
satisfies (\ref{condstar}). Under $\dagger $ the
RHS of relations (\ref{def4.1}), (\ref{def4.2}) are mapped
into the RHS of relations (\ref{def4.3}), (\ref{def4.4}),
provided 
\eq
(\aa^i)^{\dagger }=\aa^+_i;
\label{star3}
\en
in this case we find, as requested
\eq
(A^i)^{\dagger }= A^+_i.
\label{star2}
\en
If $\g=sl(N),so(N)$ and $\rho=\rho_d$ 
condition (\ref{star3}) is satisfied by choosing
\eq
v^{-1}=u=\cases{\sqrt{y_{sl(N)}}
\z \mbox{if \g=sl(n)}\cr
\sqrt{y_{so(N)}}
\z \mbox{if \g=so(n)}\cr}
\en

${\cal A}_{+,so(N),\rho_d}$ admits also an alternative
$*$-structure compatible with $\tr$, namely
$(\tilde A^+_i)^{\dagger }=\tilde  A^+_jC^{ji}$
\cite{frt}
together with a nonlinear equation for
$(\tilde A^i)^{\dagger_h}$ \cite{oleg2}
which we omit here;
in this case one usually denotes the generators
by $X_i,\partial^i$ instead of $\tilde A^+_i ,\tilde A^i$,
because in the classical limit they become the
Cartesian coordinates and partial derivatives
of the $N$-dim Euclidean space respectively.
The classical limit of this $\dagger_h$ is
\eq
(a^+_i)^{\dagger }= a^+_jc_{ji} \z\z
(a^i)^{\dagger }= - c_{ij}a^j;
\label{star4}
\en
using relations (\ref{star4}), (\ref{partico}), 
$tr(\rho_d)=0$, $S x=-x$ if $x\x\in g$, one finds
again relations (\ref{starhom}).
$\tr$ as defined in formula  (\ref{defprop})  also
satisfies (\ref{condstar}). Under $\dagger $ the
RHS of relation (\ref{def4.1}) is mapped
into the RHS of relations (\ref{def4.2}), 
provided that $(\aa^+_i)^{\dagger }= \aa^+_jc_{ji}$, 
\ie
\eq
v=1\z\z u=y_{so(N)};
\en
in this case we find, as requested
\eq
(A^+_i)^{\dagger }= A^+_jC_{ji},
\en
and it is not difficult to show that 
$(A^i)^{\dagger }$ is the (nonlinear) function of
$A^i,A^+_i$ which was found in Ref. \cite{oleg2}.

\renewcommand{\theequation}{\thesection.\arabic{equation}}
\sect{Outlook, final remarks and conclusions}
\label{conclu}

Given some solutions $\aa^i,\aa^+_i$
[in the form (\ref{defaa})] of equations 
(\ref{figata1}-\ref{figata3}), the 
$A^i,A^+_i$ defined through formulae 
(\ref{def3}) (where we choose a minimal $\F$) satisfy 
the deformed commutation relations of $\Aq$ and are covariant
under the $\uqg$ action $\tr$ defined in formula
(\ref{defprop}). The corresponding basic 
algebra homomorphism $f: \Aq\rightarrow \A[[h]]$ is
defined iteratively starting from 
$f(\tilde A^i):=A^i$, $f(\tilde A^+_i):=A^+_i$.
Explicit solutions $\aa^i,\aa^+_i$ of equations 
(\ref{figata1}-\ref{figata3}) are given by
\begin{itemize}
\item formulae (\ref{ansatz0}), (\ref{yslN})
for ${\cal A}^h_{\pm,sl(N),\rho_d}$;
\item formulae (\ref{ansatz}),
(\ref{solve}) for ${\cal A}^h_{+,so(N),\rho_d}$.
\end{itemize}
The main feature of such a realization $(f,\tr)$ is that 
the $\g$-invariant ground state $|0\rangle$ as well as the first excited
states $a^+_i|0\rangle$ of the classical Fock space representation
are also respectively $U_h\g$-invariant
ground state $|0_h\rangle$ and first excited
states $A^+_i|0_h\rangle$ of the deformed Fock space representation.

According to relation
(\ref{innerauto}), all other elements of
$\A[[h]]$ 
satisfying the DCR of $\Aq$ can be
written in the form
\eq
A^{\alpha\,i}=\alpha\,A^i\alpha^{-1}\qquad\qquad
A^{+}_{\alpha\,i}= \alpha\, A^+_i\, \alpha^{-1},
\label{aa'}
\en
with $\alpha=\idA+O(h)\in\A[[h]]$.
They are manifestly covariant under the \uqg-action
$\tr_{\alpha}$ defined by
\eq
x\,\trc_{h,\alpha}\,a\: :=\:\alpha\sigma_{\varphi_h}(x_{(\bar 1)}) 
\, a\,\sigma_{\varphi_h}(x_{(\bar 2)})\alpha^{-1}.
\label{tr'}
\en
For these realizations the deformed ground state in the Fock space 
representation
reads $|0_h\rangle=\alpha|0\rangle$; thus in general 
the $\g$-invariant ground state and first excited
states of the classical Fock space representation
do not coincide with their deformed counterparts.

In this way we have found all possible pairs 
$(f_{\alpha},\trc_{h,\alpha})$ making
the diagram (\ref{diagram}) in the introduction commutative.

Note that the change
$\varphi_h\rightarrow\varphi_{h,v}=v\,\varphi_h(\cdot)\,v^{-1}$
[formula (\ref{simtra})] of the algebra isomorphism
$\uqg\rightarrow U\g[[h]]$ 
 amounts to the particular transformation 
$(f,\tr)\rightarrow (f_{\alpha},\trc_{h,\alpha})$, 
with $\alpha=\sigma(v)$.

In Sect \ref{inva} we have shown {formula (\ref{def9})]
how to construct \g-invariants $I_h^n\in\A[[h]]$ in the form of
homogeneous polynomials in $A^i,A^+_i$. 
It is immediate to verify that under a  
transformation $(f,\tr)\rightarrow (f_{\alpha},\trc_{h,\alpha})$
these $I_h^n$  transform into
$I_h^{\alpha\,n}:=\alpha\,I_h^n\,\alpha^{-1}$.

In Sect. \ref{star} we have shown (sticking to 
the explicit case of ${\cal A}^h_{\pm,sl(N),\rho_d}$ and
${\cal A}^h_{+,so(N),\rho_d}$) that, if $\Aq$ is a
module $*$-algebra [formula (\ref{condstar})]
of the compact section of \uqg ($q>1$), then 
one can choose $(f,\tr)$ so that $f$ is a $*$-homorphism,
$f(b^{\dagger_h})=[f(b)]^{\dagger }$, and $\tr$
also satisfies equation (\ref{condstar}).
It is straighforward to verify that $(f_{\alpha},\trc_{h,\alpha})$
satisfy the same constraints
provided that $\alpha$ is ``unitary'',
$\alpha^{\dagger}=\alpha^{-1}$.

\medskip

Summing up, in the present 
work we have shown how to realize a deformed $U_h\g$-covariant
Weyl or Clifford algebra $\Aq$ within the undeformed one $\A[[h]]$.
Given  a deforming map $f$ and a representation $(\pi,V)$
of $\A$ on a vector space $V$, does $(\pi\circ f,V)$ 
provide also a representation
of $\Aq$? In other words, can one interpret the elements of $\Aq$ 
as operators acting on $V$, if the elements of $\A$ are? If so,
which specific role play the elements $A^i,A^+_i$ of $\A[[h]]$?

In view of the specific example we have examined in ref. \cite{fio}
the answer to the first question seems to be always positive, whereas
the
converse statement is wrong: 
\eg for any $h\in\b{R}$ there are more (inequivalent)
representations of the deformed algebra than representations of
the undeformed one. This may seem a paradox, because 
in a $h$-formal-power sense $f^{-1}$ can be defined and gives
$f^{-1}(a^i)=\tilde A^i+O(h)$, $f^{-1}(a^+_i)=\tilde A^+_i+O(h)$,
whence at least for small $h$ one would expect the 
deformed and undeformed representation 
theories to coincide; but in fact there is no warrancy that,  
also in some operatorial sense, for small $h$ the would-be
$f^{-1}(a^i),f^{-1}(a^+_i)$ are `close' to $\tilde A^i,\tilde A^+_i$.
Of course, we are especially interested in Hilbert
space representations of $*$-algebras:
in Ref. \cite{fio} we checked that in the operator-norm topology
$f^{-1}$ is ill-defined on all but exactly one deformed representation.
Roughly speaking, the reason is that the `particle
number' observables $n^i$, which enter the transformation $f$ (see
e.g. (\ref{lastb})) are unbounded operators, therefore even for very 
small $h$ the effect of the transformation on their large-eigenvalue 
eigenvectors can be so large to `push' the latter out of the domain of
definition of the operators in $f^{-1}(\A)$.

We are especially interested in the case of $*$-algebras admitting Fock 
space representations. The results presented in the previous paragraphs
could in principle be applied to models in quantum field
theory or condensed matter physics by choosing representations
$\rho$ which are the direct sum of many copies of the same
fundamental representation $\rho_d$; this is what we have 
addressed in Ref. \cite{fionuovo}. The different copies would 
correspond respectively to different space(time)-points or
crystal sites. 

One important physical issue is if $U_h\g$-covariance necessarily
implies exotic particle statistics. In view of what we have said
the answer is no \cite{fiogrou}. 
At least for compact $\g$ and $U_h\g$ ($h$ is real),
the undeformed Fock space representation, which allows a
`Bosons \& Fermions' particle interpretation, would carry also a
representation of the deformed one. Next point is the role
of the operators $A^i,A^+_j$. Quadratic commutation`r%lations
of the type (\ref{gqcr0}-\ref{gqcr2}) mean  
that $A^+_i,A^i$ act as creators and
annihilators of some excitations; a glance at (\ref{def3}), (\ref{aa'})
shows that these are not the undeformed excitations,
but some `composite' ones. The last point is:
what could these operators be good for. As an Hamiltonian $H$ of the
system we may choose a simple combination of the $U_h\g$-invariants
$I^n_h$ of section \ref{inva}; thus the Hamiltonian is $U_h\g$-invariant 
and has a simple polynomial structure in the composite operators
$A^i,A^+_j$. $H$ is also $\g$-invariant, but has a highly non-polynomial
structure in the undeformed generators $a^i,a^+_j$ (it would be
tempting to understand what kind of physics it could describe!).
This suggests that the use of the $A^i,A^+_j$ instead of the $a^i,a^+_j$
should simplify the resolution of the corresponding dynamics
(similarly to what has been suggested in Ref. \cite{wecorfu} for
a 1-dim toy-model).

\renewcommand{\theequation}{\thesubsection.\arabic{equation}}
\app{Appendix}

\subsect{Proof of Lemma \ref{lemma2}}

\begin{lemma}\cite{fio}
If ${\cal T}\in U\g[[h]]^{\ot^3}$ 
is \g-invariant (\ie $[{\cal T},U\g[[h]]^{\ot^3}]=0$) then 
$m_{ij}S_{i}{\cal T}$, $m_{ij}S_{j}{\cal T}$ ($i,j=1,2,3$, $i\neq j$)
are $\g$-invariants belonging to $U\g[[h]]^{\ot^2}$.
\label{lemmino}
\end{lemma}
\vskip-.05cm
 (Here $S_{i}$ denotes $S $ acting on
the $i$-th tensor factor, and 
$m_{ij}$ multiplication of the $i$-th tensor factor by
the $j$-th from the right.)

We may apply the previous lemma to
 ${\cal T}=\phi$, or ${\cal T}=\phi^{-1}$. Looking at 
the definition (\ref{defphi}) one finds in particular the
following $\g$-invariants
\eq
\begin{array}{rcl}
T_1 &:=& m_{12}S_{1}\phi=(S \F^{(1)}\gamma\ot \1)\F
(\F^{(2)}_{(1)}\ot \F^{(2)}_{(2)}),\\
T_2 &:=& m_{23}S_{3}\phi=
(\F^{-1(1)}_{~(1)}\ot \F^{-1(1)}_{~(2)})\F^{-1}(\1\ot\gamma^{-1}
S \F^{-1(2)}) ;
\end{array}
\label{quanti!}
\en
alternative expressions for these $T_i$ can be obtained by applying
the same operations to the identities
\eqa
\phi  q^{t_{12}+t_{13}\over 2}& = & q^{t_{12}\over 2} \F^{-1}_{23,1}
\F^{-1}_{12}\R_{13}\F_{32}\F_{1,23},\label{uno}\\
q^{-\frac{t_{13}+t_{23}}2}\phi & = &
\F^{-1}_{12,3}\F^{-1}_{21}\R^{-1}_{13}\F_{23}\F_{3,12}
q^{-\frac{t_{23}}2},
\label{due}
\ena
which directly follow from relations 
(\ref{simpler}), (\ref{defphi}),
(\ref{defR}) and the observation
that $[\phi,q^{t_{12}+t_{13}+t_{23}\over 2}]=0$. 
Applying $m_{12}S_{1}$ to (\ref{due}), 
$m_{23}S_{3}$ to (\ref{uno}) we get
\eq
\begin{array}{rcl}
T_1 &=& (S \F^{(2)}\gamma'^{-1}\ot \1)\F_{21}
(\F^{(1)}_{(1)}\ot \F^{(1)}_{(2)}),\\
T_2 &=& 
(\F^{-1(2)}_{~(1)}\ot \F^{-1(2)}_{~(2)})\F^{-1}_{21}(\1\ot\gamma'
S \F^{-1(1)}).
\end{array}
\label{tanti!}
\en
>From eq. (\ref{quanti!}), (\ref{tanti!}) 
we easily find out that the inverse of $T_i$ take the form
\eqa
T_1^{-1} &=& \F^{-1}\left[\gamma^{-1}(S \F^{-1(1)})
\F^{-1(2)}_{~(1)}\ot \F^{-1(2)}_{~(2)}\right]\label{tantancora!1}\\
 & = &  \F^{-1}_{21}\left[\gamma'(S \F^{-1(2)})
\F^{-1(1)}_{~(1)}\ot \F^{-1(1)}_{~(2)}\right],\label{tantancora!2}\\
T_2^{-1} &=& \left[\F^{(1)}_{(1)}\ot \F^{(1)}_{(2)}
(S \F^{(2)})\gamma\right]\F \label{tantancora!3}\\
 &= & \left[\F^{(2)}_{(1)}\ot \F^{(2)}_{(2)}
(S \F^{(1)})\gamma'^{-1}\right]\F_{21},
\label{tantancora!4}
\ena
since 
$[T_i,\F^{\pm 1 (j)}_{(1)}\ot \F^{\pm 1 (j)}_{(2)}]=0$,
with $i,j=1,2$.

If $\F$ is minimal 
then it is easy to verify that, according to (\ref{lilla}) and their 
definitions (\ref{quanti!}), $T_i\equiv\1^{\ot^3}$. In the latter case
the last four relations are equivalent to relations 
(\ref{dritto1}-\ref{dritto4}).

\subsect{Proof of Proposition \ref{propinv}}

We start by expressing 
$(A...A)^{I_n}\equiv A^{i_1}...A^{i_{h_n}}$,
$(A^+...A^+)_{J_n}\equiv A^+_{j_1}...A^+_{j_{k_n}}$ respectively
in the form $(\aa...\aa)\sigma(\cdot)$, $(\aa^+...\aa^+)\sigma(\cdot)$.
First note that 
\eqa
A^{i_1}A^{i_2} &\stackrel{(\ref{def4.4})}{=}&
\rho(\gamma^{-1}S\F^{-1(1)})^{i_1}_{l_1}
\rho(\gamma^{-1}S\F^{-1(1')})^{i_2}_{l_2}
\aa^{l_1}\sigma(\F^{-1(2)})\aa^{l_2}\sigma(\F^{-1(2')})
\nonumber
\\
& = &\rho^\vee(\F^{-1(1)}\gamma')_{i_1}^{l_1}
\rho^\vee(\F^{-1(1')}\gamma')_{i_2}^{l_2}
\aa^{l_1}\sigma(\F^{-1(2)})\aa^{l_2}\sigma(\F^{-1(2')})
\nonumber\\
& \stackrel{(\ref{cov}),(\ref{covl})}{=}&
\rho^\vee(\F^{-1(1)}\gamma')_{i_1}^{l_1}
\rho^\vee(\F^{-1(2)}_{(1)}\F^{-1(1')}\gamma')_{i_2}^{l_2}
\aa^{l_1}\aa^{l_2}\sigma(\F^{-1(2)}_{(2)}\F^{-1(2')})
\nonumber\\
& = &
\aa^{l_1}\aa^{l_2}\left[
\left(\rho^{\vee l_1}_{~i_1}\ot \rho^{\vee l_2}_{~i_2}
\ot \sigma\right)\left(\F^{-1}_{1,23}\F^{-1}_{23}
(\gamma'^{\ot^2}\ot\1)\right)\right]
\nonumber
\ena
whence, by repeated application, we find
\eqa
(A...A)^{I_n} &\stackrel{(\ref{def10})}{=}&
(\aa...\aa)^{L_n}\left[
\left((\rho^{\vee\ot^{h_n}})^{L_n}_{I_n}\ot\sigma\right)
\left(\F^{-1}_{12...(h_n\!+\!1)}
(\gamma'^{\ot^{h_n}}\ot\1)\right)\right];\nonumber
\ena
similarly, starting from relation (\ref{def4.1}) we find
\eqa
(A^+...A^+)_{J_n} & = &
(\aa^+...\aa^+)_{M_n}\left[
\left((\rho^{\ot^{k_n}})^{M_n}_{J_n}\ot\sigma\right)
\left(\F^{-1}_{12...(k_n\!+\!1)}\right)\right].\nonumber
\ena
Putting these results together we find
\[
(A^+...A^+)_{J_n}(A...A)^{I_n}\stackrel{(\ref{covl}),
(\ref{cov})}{=} \z\z\z\z\z\z\z\z\z\x\x
\]
\[
(\aa^+...\aa^+)_{M_n}(\aa...\aa)^{L_n}
\left[\left((\rho^{\ot^{k_n}})^{M_n}_{J_n}\ot
(\rho^{\vee\ot^{h_n}})^{L_n}_{I_n}\ot\sigma\right)
\left(\F^{-1}_{12...(b_n\!+\!1)}
(\1^{\ot^{k_n}}\ot\gamma'^{\ot^{h_n}}\ot\1)\right)\right],
\]
whence,
\eq
I_h^n  \stackrel{(\ref{propto})}{\propto} 
(\aa^+...\aa^+)_{M_n}(\aa...\aa)^{L_n}
\left[\left((\rho^{\ot^{k_n}})^{M_n}_{J_n}\ot
(\rho^{\vee\ot^{h_n}})^{L_n}_{I_n}\ot\sigma\right)
\F^{-1}_{12...(b_n\!+\!1)}\F_{12...b_n,b_n\!+\!1}
\right]d^{L_n}_{M_n}.
\label{pussi}
\en
We prove now that 
\eq
\F^{-1}_{12...(b\!+\!1)}\F_{12...b}=
\phi^{-1}_{(b\!-\!1)b(b\!+\!1)}
\phi^{-1}_{(b\!-\!2),(b\!-\!1)b,(b\!+\!1)}
...\phi^{-1}_{1,2...b,(b\!+\!1)}\F^{-1}_{12...m,m\!+\!1};
\label{billy}
\en
then the claim will follow from relation (\ref{pussi})
and the observation that
\eq
\left[\left(\rho^{\ot^{k_n}}
\!\ot \!\rho^{\vee~\ot^{h_n}}\!\ot \!\id\right)
(\F_{12...b_n,b_n\!+\!1})
\right]^{ M_n L_n}_{ J_n I_n}d^{L_n}_{M_n}
\stackrel{(\ref{cisotrop})}{=}
\varepsilon(\F^{(1)})\F^{(2)}d^{I_n}_{J_n}
\stackrel{(\ref{cond2})}=d^{I_n}_{J_n}.
\en
To prove relation (\ref{billy}) we start from
\[
\phi^{-1}_{123}\F^{-1}_{12,3}\stackrel{(\ref{defphi})}=
\F^{-1}_{1,23}\F^{-1}_{23}\F_{12}\stackrel{(\ref{def10})}{=}
\F^{-1}_{123}\F_{12}; 
\]
this is relation (\ref{billy}) for $b=2$.
Applying $\id\ot\Delta\ot \id$ and multiplying 
the result from the left by $\phi^{-1}_{234}$  we find
\eqa
\phi^{-1}_{234}\phi^{-1}_{1,23,4}\F^{-1}_{123,4} & = &
\phi^{-1}_{234}\F^{-1}_{1,234}\F^{-1}_{2,34}\F_{1,23}\nonumber\\
&\stackrel{(\ref{ginv})}{=}&\F^{-1}_{1,234}
\phi^{-1}_{234}\F^{-1}_{2,34}\F_{1,23}\nonumber\\
&\stackrel{(\ref{defphi})}{=}&
\F^{-1}_{1,234}\F^{-1}_{2,34}\F^{-1}_{34}\F_{23}\F_{1,23},\nonumber\\
&\stackrel{(\ref{def10})}{=}&
\F^{-1}_{1233}\F_{123},
\nonumber
\ena
\ie relation (\ref{billy}) for $b=3$. Applying 
to the latter relation $\id\ot\Delta^{(2)} \ot \id$ and multiplying 
the result from the left by $\phi^{-1}_{345}$ we find 
relation (\ref{billy}) for $b=4$, and so on $\Box$.

\subsect{Proof of Proposition \ref{figate}}

\eqa
0&\stackrel{(\ref{gqcr0})}{=} 
& A^i A^j \mp \P^{ji}_{hk} A^kA^h  
\nonumber \\
 &\stackrel{(\ref{def4.3})}{=} & (\1\mp \P)^{ji}_{hk}
\rho(\F^{(1)})^h_l\rho(\F^{(1')})^k_m\sigma(\F^{(2)}) \aa^m
\sigma(\F^{(2')})\aa^l
\nonumber \\
& \stackrel{(\ref{cov}),(\ref{covl})}{=}&(\1\mp \P)^{ji}_{hk}
\rho(\F^{(1)}\F^{(2')}_{(2')})^k_m\rho(\F^{(1)})^h_l
\sigma(\F^{(2)}_{(1)}\F^{(2')}_{(1')})\aa^m \aa^l;
\nonumber 
\ena
multiplying both sides from the left by 
$(\rho\ot\rho\ot\sigma)(\F_{1,23}^{-1}\F_{23}^{-1})$
and noting that
\[
\P_{12}\stackrel{(\ref{genperm}), (\ref{UVmatrix})}{=}
[(\rho\ot\rho\ot\sigma)\F_{12}\F_{12,3}]U_{12}
[(\rho\ot\rho\ot\sigma)\F_{12,3}^{-1}\F_{12}^{-1}],
\]
we find
\\
\[
\{\1^{\ot^3}\mp [(\rho^{\ot^2}\ot \sigma)
\F_{1,23}^{-1}\F_{23}^{-1}\F_{12}\F_{12,3}]
U_{12}[(\rho^{\ot^2}\ot \sigma)\F_{12,3}^{-1}\F^{-1}_{12}\F_{23}
\F_{1,23})]\}^{hk}_{lm}\aa^m \aa^l=0,
\]
\ie relation (\ref{figata1}), once we take
definitions (\ref{def8}), (\ref{defphi}) into account.
Using definition (\ref{def4.1}) one can prove
in a similar way that relations (\ref{gqcr1}),
(\ref{figata2}) are equivalent. Similarly, 
\eqa
0&\stackrel{(\ref{gqcr2})}{=} 
& A^i A^+_j -\delta^i_j\idA
\mp \tilde{P}^{F~ih}_{~~~jk} A^+_h A^k 
\nonumber \\
 &\stackrel{(\ref{def4.3}),(\ref{def4.1})}{=} & 
\rho(\F^{(1)})^i_l\sigma(\F^{(2)})\aa^l\aa^+_m
\sigma(\F^{-1(2')})\rho(\F^{-1(1')})^m_j 
\nonumber \\
& & -\delta^i_j\idA
\mp \aa^+_m \sigma(\F^{-1(2)})\rho(\F^{-1(1)})^m_h     
\tilde{P}^{F~ih}_{~~~jk}\rho(\F^{(1')})^k_l\sigma(\F^{(2)})\aa^l; 
\nonumber 
\ena
multiplying both sides by $\rho(\F^{-1(1)})^{i'}_i\sigma(\F^{-1(2)})$
from the left and by $\rho(\F^{(1')})^j_{j'}\sigma(\F^{(2')})$
from the right, and noting that
\[
\tilde{P}_{12}^F\stackrel{(\ref{genrmat}), 
(\ref{UVmatrix})}{=}[(\rho^{\ot^2}\ot\sigma)\F_{12}\F_{12,3}]V_{12}
[(\rho^{\ot^2}\ot\sigma)\F_{12,3}^{-1}\F_{12}^{-1}],
\]
we get
\eqa
0 & = & \aa^{i'}\aa^+_{j'}-\delta^{i'}_{j'}\idA 
\mp \rho(\F^{-1(1)})^{i'}_i\sigma(\F^{-1(2)})
\aa^+_m \left\{[(\rho^{\ot^2}\ot\sigma)\F^{-1}_{13}\F_{12}\right.
\F_{12,3}]
\nonumber \\
&&\left.\times V_{12}[(\rho^{\ot^2}\ot\sigma)\F_{12,3}^{-1}
\F_{12}^{-1}\F_{13}]
\right\}^{im}_{jl}\aa^l\rho(\F^{(1')})^j_{j'}\sigma(\F^{(2')})
\nonumber \\
& \stackrel{(\ref{cov}),(\ref{covl})}{=}&
\aa^{i'}\aa^+_{j'}-\delta^{i'}_{j'}\idA
\mp 
\aa^+_m \left\{[(\rho^{\ot^2}\ot\sigma)\F^{-1}_{1,23}
\F^{-1}_{13}\F_{12}\F_{12,3}]\right.
\nonumber \\
&&\left. \times  V_{12}[(\rho^{\ot^2}\ot\sigma)
\F_{12,3}^{-1}\F_{12}^{-1}\F_{13}F_{1,23}]\right\}^{i'm}_{j'l}\aa^l,
\nonumber 
\ena
whence the equivalence between relations (\ref{gqcr2}),
(\ref{figata3}) follows, once one recalls the definition
(\ref{defphi}). $\Box$ 

\subsect{Some properties of special and q-special functions}

The following results can be found in standard textbooks.
If the parameters $a,b,c\in\cn$ are such that none
of the quantities $c\!-\!1,a\!-\!b,a\!+\!b\!-\!c$ 
is a positive integer, the general solution of the
hypergeometric differential equation 
in the complex $z$-plane 
\eq
y''(1-z)z+y'[c-(a\!+\!b\!+\!1)z]-yab=0
\label{hyper}
\en
can be expressed as some combinations
\eqa
&&\!\!y(z) = \alpha\, F\!(a,b,c; z)+\beta\, z^{1\!-\!c}
F(1\!+\!a\!-\!c,1\!+\!b\!-\!c,2\!-\!c; z),
\label{combi1}\\
&&\!\!=\gamma\, F\!(a,b,a\!+\!b\!+\!1\!-\! c;1\!-\!z)
 +\delta\,
(1\!-\!z)^{c\!-\!a\!-\!b}\!
F(c\!-\!a,c\!-\!b,c\!+\!1\!-\!a\!-\!b;1\!-\!z),\z\x\:
\label{combi2}
\ena
where $\alpha,\beta, \gamma,\delta\in\cn$ and
$F(a,b,c; z)$ is the hypergeometric function.
As known, 
\eqa
F(a,b,c; 0) & = & 1,             \label{value} \\
{d\over dz}F(a,b,c; z) & = & {a\, b\over c}
F(a\!+\!1,b\!+\!1,c\!+\!1; z). \label{deri}
\ena
The combinations (\ref{combi1}), (\ref{combi2}) explicitly
display the singular and non-singular part of the solution
respectively around the poles $x=0,1$. 
An essential identity to determine the asymptotic behaviour of a
solution $y$ around the pole $x=0$ (resp. $x=1$), known its
asymptotic behaviour around the pole $x=1$ (resp. $x=0$), is
\eqa
F(a,b,c; z) \!&\! =\! &\! {B(c,c\!-\!a\!-\!b)\over 
B(c\!-\!a,c\!-\!b)}
F(a,b,a\!+\!b\!+\!1\!-\!c;1\!-\! z)
\label{lintra} \\
\!&\!+\!&\! {B(c,a\!+\!b\!-\!c)\over B(a,b)}
(1\!-\!z)^{c\!-\!a\!-\!b}F(c\!-\!a,c\!-\!b,c\!+\!1\!-\!a\!-\!b;
1\!-\!z);\nonumber
\ena
Here $\Gamma(a)$ and
$B(a,b)$ are Euler's $\Gamma$- and $\beta$-functions respectively;
as known, 
\eq
\Gamma(a\!+\!1) = a\Gamma(a)  \z\z  
B(a,b) ={\Gamma(a)\Gamma(b)\over \Gamma(a\!+\!b)}.
\label{betagamma}
\en

A less obvious property is 
\eq
\Gamma(a)\Gamma(-a)=-\frac {\pi}{a\sin{\pi a}}.
\label{bingo}
\en
The $q$-gamma function $\Gamma_q$ can be defined
when $|q|<1$ by \cite{Exton}
\eq
\Gamma_q(a):= (1-q^{1-a})\prod\limits_{k=0}^{\infty}
\frac{(1-q^{k\!+\!1})}{(1-q^{a\!+\!k})}=(1-q^{1-a})
\sum\limits_{n=0}^{\infty}\frac{(q^{1-a};q)_n}{(q^a;q)_n}q^{na},
\label{defgammaq}
\en
where $(a;q)_n:=\prod\limits_{k=0}^{n-1}(1-aq^k)$;
it satisfies the following modified version of the property
(\ref{betagamma})$_1$:
\eq
\Gamma_q(a\!+\!1)= (a)_q\Gamma_q(a), \z\z\z 
(a)_q:={(q^a-1)\over (q-1)}.
\label{gammaq}
\en
We introduce also a different version of the $q$-gamma function
by 
\eq
\tilde \Gamma_q(a):=\Gamma_{q^2}(a)q^{-\frac{a(a-3)}2};
\label{defgammaq'}
\en
the latter satisfies
\eq
\tilde \Gamma_q(a\!+\!1)= [a]_q\tilde \Gamma_q(a), \z\z\z 
[a]_q:={(q^a\!-\!q^{-a})\over (q\!-\!q^{-1})}.
\label{gammaq'}
\en

\subsect{Proof of Theorem \ref{propsln}}

$Proof$.
We need to show that equations (\ref{figata1}-\ref{figata3})
are fulfilled. We get rid of
indices by introducing the following vector notation:
\eqa
&(aa)^{ij} := a^ia^j  \z\z &(a^+a^+)_{ij} := a^+_ia^+_j \nonumber\\
&(av)^{ij} := a^iv^j \z\z& (va)^{ij} := v^i a^j \nonumber\\
&(a^+w)_{ij} := a^+_iw_j \z\z& (wa^+)_{ij}:= w_ia^+_j\nonumber\\
&w\cdot a := w_ia^i \z\z& a^+\cdot v := a^+_iv^i,\nonumber
\ena
where $v\equiv (v^i)\in \cn^N$, $w\equiv (w_i)\in \cn^N$
denote arbitrary covariant and controvariant vectors 
respectively.  If we plug  (\ref{def13}) into
(\ref{figata1}-\ref{figata3}), factor out of 
(\ref{figata1}) and (\ref{figata2}) $I(n)\,I(n\!+\!1)$
and $\tilde I(n)\,\tilde I(n\!+\!1)$ respectively, 
multiply eq. (\ref{figata3}) by $v^jw_i$, then we find the
equivalent system (in vector notation)
\eqa
aa\, & = & \:\pm\: (M^{-1}P M)aa  \label{figat1}\\
a^+a^+\,& = & \:\pm\: a^+a^+\,(M^{-1}PM)
\label{figat2}\\
{(n\!+\!1)_{q^{\pm 2}}\over (n\!+\!1)}
(w\!\cdot\! a)\cdot(a^+\!\cdot\! v) & = & w\!\cdot\! v\,\idA\:\pm\: 
q^{\pm 1}{n_{q^{\pm 2}}\over n}\,wa^+\,(M^{-1}\,V\,M)\,va.
\label{figat3}
\ena
It is straightforward to show that
\eqa
&A\,va  =  \:\pm\:(n\!-\!1)\, va  \z &
a^+\! w\, A  =  -wa^+\mp a^+\!w +(w\!\cdot\!a)\,a^+\!a^+ \nonumber\\ 
&  wa^+\,P  =  a^+\! w \z & a^+\! w \,P  =  wa^+. \label{utile1} 
\ena
As a consequence one finds, in particular,
\eqa
a^+\! w\, (A\!+\!P)^k &= & (\mp 1)^k a^+\! w\,\pm
\frac{(\pm n)^k-(\mp 1)^k}{n+1}\,
(w\!\cdot\!a)\,a^+\!a^+\z\Rightarrow
\nonumber\\
\Rightarrow \z a^+\! w\,q^{A\!+\!P} &=& q^{\pm 1}
a^+\! w\,\pm
\frac{q^{\mp n}-q^{\mp 1}}{n+1}\,(w\!\cdot\!a)\,a^+\!a^+
\label{utile2}
\ena

Let us prove eq. (\ref{figat1}), (\ref{figat2}). The matrix
$M$ (\ref{def8}) takes the form
\eq
M \stackrel{(\ref{integral})}{=}
\lim_{x_0,y_0\rightarrow 0^+}\left\{x_0^{-2\hbar P}
\vec{P}\exp\left[-2\hbar\int\limits^{1-y_0}_{x_0}dx
\left({P\over x}+
{A\over x-1}\right)\right] y_0^{2\hbar A}\right\};
\label{gigg}
\en
the contributions  of the central terms $-2\frac {\1^{\ot^3}}N$,
$-\frac2N \1^{\ot^3}\ot n$
to the integral are cancelled by the corresponding 
contributions from $x_0^{-\hbar t_{12}}$,  $y_0^{\hbar t_{23}}$,
in the limit $x_0,y_0\rightarrow 0^+$. Since $aa$ is an `eigenvector'
both of $A$ and $P$, the path-order $\vec{P}$ becomes redundant and we 
find that $M$ acts trivially on $M$:
\eqa
Maa & = & aa\,
\lim_{x_0,y_0\rightarrow 0^+}\left\{x_0^{\mp2\hbar}
\exp\left[-2\hbar\int\limits^{1-y_0}_{x_0}dx
\left(\pm {1\over x}\pm
{n-2\over x-1}\right)\right] y_0^{\pm2\hbar(n\!-\!2)}\right\}\nonumber\\
& = & aa\,\lim_{x_0,y_0\rightarrow 0^+}(1-y_0)^{\mp 2\hbar}
(1-x_0)^{\pm 2\hbar(n\!-\!2)}=aa.
\ena
Therefore $MPM\,aa=\pm\:aa$, Q.E.D. Similarly one proves eq. 
(\ref{figat2}).

In order to prove eq. (\ref{figat3}) it is convenient to recast 
$M^{-1}\, V\, M$ in a more manageable form. Permuting the 
second and third tensor factor in eq. $(\ref{simpler})_{(1)}$, we find 
\eq
\phi^{-1}_{213}q^{t_{12}\over 2}\phi_{123}=
q^{t_{12}+t_{13}\over 2} \phi_{132}q^{- \frac{t_{23}}2},
\label{manage}
\en
whence 
\eq
M^{-1}\, V\, M=PM_{21}^{-1}q^{\pm 1+\frac 1N+\rho_d{}^{\ot^2}(\frac t2)} 
M\stackrel{(\ref{lilla2}),(\ref{def14})}{=}
Pq^{A+P}\lim_{x\rightarrow 0^+}x^{-2B}f(x)q^{-A}, 
\label{stucco}
\en
where $f$ is the $\sigma(Usl(N))[[h]]$-valued $N^2\times N^2$ matrix
satisfying the differential equation and asymptotic conditions
\eq
f'=2\hbar\left(\frac Bx+\frac A{x-1}\right)f \z\z
\lim_{x\rightarrow 1} f(x)(1-x)^{-2A}=1
\label{sskzeq}
\en
[the latter are obtained from eq. (\ref{skzeq}) by permuting the 
second and third tensor factor and by
getting rid of the central terms involved in
$(\rho\ot\rho\ot\sigma) (t_{ij})$ (formulae
(\ref{def14})) since, as  in formula (\ref{gigg}),
the latter cancel with each other in the limit $x\rightarrow 0$].

It is convenient to introduce in $\A[[h]]$ a grading $g$, by
setting $g(b)=l\in\zn$ iff $[n,b]=lb$, $b\in\A[[h]]$.
Since $g(fva)=-1$, and $fv^ia^j$ is a doubly contravariant tensor,
its most general expansion is
\eq
f(x)va=av f_1(x)+vaf_2(x)+aa(a^+\!\cdot\!v)f_3(x),
\label{decomp}
\en
where $f_i$ are invariants with $g(f_i)=0$; therefore
$f_i=f_i(n)$. Thus we find
\eqa
&&wa^+\cdot(M^{-1}\, V\, M)\,va  
\nonumber \\
&&\stackrel{(\ref{utile1}),(\ref{stucco})}{=}
a^+\! w\,q^{A+P}\lim_{x\rightarrow 0^+}x^{- 2\hbar B}
f(x)q^{-A}\,va  
\nonumber \\
&&\stackrel{(\ref{utile1}),(\ref{utile2})}{=}
\lim_{x\rightarrow 0^+}x^{\mp 2\hbar(n\!-\!1)}
\left[q^{\mp 1} a^+\! w\,\pm
\frac{q^{\pm n}-q^{\mp 1}}{n+1}\,(w\!\cdot\!a)\,a^+\!a^+
\right]\cdot(f(x)\,va)\,q^{\mp(n\!-\!1)} 
\nonumber \\
&&\stackrel{(\ref{decomp})}{=}q^{\mp(n\!-\!1)}
\lim_{x\rightarrow 0^+}x^{\mp 2\hbar(n\!-\!1)}
\left[\!q^{\mp 1}\! a^+\! w\pm
\frac{q^{\pm n}\!-\!q^{\mp 1}}{n+1}(w\!\cdot\!a)a^+\!a^+
\right]\times
\nonumber \\
&& \z\z\z\z\z\x\left[av\! f_1(x)\!
+\!va\!f_2(x)\!+\!aa(a^+\!\!\cdot\!v)\!f_3(x)\right] 
\nonumber \\
&&=q^{\mp(n\!-\!1)} 
\lim_{x\rightarrow 0^+}x^{\mp 2\hbar(n\!-\!1)}
\left\{q^{\mp 1} (w\!\cdot\!v)(nf_1\mp f_2)+
(w\!\cdot\!a)(a^+\!\!\cdot\! v)\left[q^{\mp 1}
(nf_3\pm f_2)\right.\right. 
 \nonumber \\
&&\left.\left. \z\z\z\z\z\x +  n\frac{q^{\pm n}-q^{\mp 1}}{n+1}
\left(f_1\pm f_2+ (n\!+\!1)f_3\right)\right]\right\}
\nonumber \\
&&=q^{\mp n} \left\{(w\!\cdot\!v)\,
l_1+(w\!\cdot\!a)(a^+\!\!\cdot\! v)\left[l_2
+  n\frac{q^{\pm (n\!+\!1)}-1}{n+1}\,l_3\right]\right\},
\label{brodo}
\ena
where we have defined
\eqa
l_1 &:=& \lim_{x\rightarrow 0^+}x^{\mp 2\hbar(n\!-\!1)}
(nf_1(x)\mp f_2(x))\nonumber \\
l_2 &:=& \lim_{x\rightarrow 0^+}x^{\mp 2\hbar(n\!-\!1)}
(nf_3(x)\pm f_2(x))\nonumber \\
l_3 &:=& \lim_{x\rightarrow 0^+}x^{\mp 2\hbar(n\!-\!1)}
\left[f_1(x)\pm f_2(x)+(n\!+\!1)f_3(x)\right].
\ena

To evaluate the limits $l_i$ let us consider 
the linear system of first order differential equations
satisfied by $f_i$.
>From (\ref{sskzeq}) we find
\eqa
f_1' &= & \hbar\left[\pm\left({1\over 1-x} +{n-1\over x}
\right)f_1-{f_2\over x}\right] \label{syst1}\\
f_2' &= & \hbar\left[{f_1\over 1-x}\mp\left({n-1\over 1-x} 
+{1\over x}\right)f_2\right] \label{syst2}\\
f_3' &= & \hbar\left[\mp {f_1\over 1-x}+{f_2\over x}
\mp\left({1\over 1-x} 
-{1\over x}\right)(n\!-\!1)f_3\right] \label{syst3}
\ena
and the asymptotic condition
\eq
\lim_{x\rightarrow 1}f_1(x)=0=\lim_{x\rightarrow 1}f_3(x)
\z\lim_{x\rightarrow 1}f_2(x)(1-x)^{\mp 2\hbar(n\!-\!1)}=1.
\label{asymp3}
\en

The first two equations can be solved separatly; 
then the third will yield $f_3$
in terms of $f_1,f_2$ just by an integration. 
Actually one of the combination we are 
interested in, $[\!f_1(x)\!\pm\! f_2(x)\!+\!(n\!+\!1)f_3(x)]$, 
satisfies a completely decoupled equation,
\[
\frac{d}{dx}[f_1(x)\!\pm\! f_2(x)\!+\!(n\!+\!1)f_3(x)]
=\pm 2\hbar(n\!-\!1)\left[{1\over x}+{1\over x-1}\right]
[f_1(x)\!\pm\! f_2(x)\!+\!(n\!+\!1)f_3(x)],
\]
which [taking  conditions (\ref{asymp3}) into
account] is easily integrated to
\eq
f_1(x)\!\pm\! f_2(x)\!+\!(n\!+\!1)f_3(x)=\pm
[x(1-x)]^{\pm 2\hbar(n-1)}.
\en
This will yield therefore $f_3$ in terms of $f_1,f_2$.
Dividing (\ref{syst1}) by $f_1$, (\ref{syst2}) by $f_2$
we find
\eqa
{f_1'\over f_1} &= & 2\hbar\left[\pm\left({1\over 1-x} +{n-1\over x}
\right)-\frac 1x\,{f_2\over f_1}\right] \label{syste1}\\
{f_2'\over f_2} &= & 2\hbar\left[{1\over (1\!-\!x)}{f_1
\over f_2}\mp
\left({n-1\over 1-x} +{1\over x}\right)\right]\label{syste2}
\ena
taking the difference of the two, one finds a Riccati
equation in the unknwon $u:={f_1\over f_2}$:
\eq
{u'\over u}=\frac{d}{dx}\ln({f_1\over f_2})=
{f_1'\over f_1}-{f_2'\over f_2}=2\hbar\left[\pm n
\left({1\over x}+{1\over 1-x}\right)-{u^{-1}\over x}
-{u\over 1-x}\right];
\en
this should be supplemented with the condition
$u\stackrel{x\rightarrow 1}{\rightarrow}0$.
To get rid
of its nonlinearity one can transform it into a (linear)
second order equation in an unknown $y(x)$
by a standard substitution, which
in this case takes the form
\eq
u={y'\over y}{(1\!-\!x)\over 2\hbar};
\en
the new equation will read
\eq
y''(1-x)x-y'(x\pm n2\hbar)+(2\hbar)^2 y=0.
\en
We recognize the hypergeometric equation [formula
(\ref{hyper})] with parameters
\eq
a=\pm 2\hbar\z\z b=\mp 2\hbar\z\z c=\mp 2 n\hbar.
\en
Its general solution can be expressed in the form
(\ref{combi2}), in terms of the hypergeometric
function $F$. Imposing the condition
$\lim_{x\rightarrow 1}{y'(1\!-\!x)\over2\hbar y}=0$
one finds that it must be $\delta=0$, implying
\eq
{f_1\over f_2}=u=\frac{1\!-\!x}{2\hbar}\frac{d}{dx}
\ln\left[F(\pm 2\hbar,
\mp 2\hbar,1\!\pm\! 2\hbar n;1\!-\!x)\right].
\en
We can now replace this result in the RHS in
eq. (\ref{syste2}):
\eq
\frac{d}{dx}\ln(f_2) = \frac{d}{dx}\ln\left[F(\pm 2\hbar,
\mp 2\hbar,1\!\pm\! 2\hbar n;1\!-\!x)\right]\mp 2\hbar
\left({n-1\over 1-x} +{1\over x}\right);
\en
taking into account the condition (\ref{asymp3}),
the latter is integrated to 
\eq
f_2(x) = x^{\mp 2\hbar}(1-x)^{\pm2\hbar(n\!-\!1)}
F(\pm 2\hbar,\mp 2\hbar,1\!\pm\! 2\hbar n;1\!-\!x).
\en
Finally, we find
\eq
f_1(x) = u(x)f_2(x)=-\frac 1{2\hbar}F'(\pm 2\hbar,\mp 
2\hbar,1\!\pm\! 2\hbar n;1\!-\!x) 
x^{\mp 2\hbar}(1-x)^{1\pm2\hbar(n\!-\!1)}.
\en

>From properties (\ref{value}),(\ref{deri}) we can 
easily read off
the asymptotic behaviour of $f_2,f_1$ for 
$x\rightarrow 0^+$:
\[
\begin{array}{l}
f_2(x)\!\!\stackrel{(\ref{lintra})}{=}
x^{\mp 2\hbar}(1-x)^{\pm2\hbar(n\!-\!1)}
\left[{B(1\!\pm\! 2\hbar n,1\!\pm\! 2\hbar n)\over
B(1\!\pm\! 2\hbar (n\!+\!1),1\!\pm\! 2\hbar (n\!-\!1))}
{\scriptstyle F\left(\!\pm\! 2\hbar,\!\mp\! 2\hbar,\!\mp\! 
2\hbar n;x\right)} \: +\: x^{1\pm 2\hbar n}\times
\right.        \nonumber\cr
\left.      
{B(1\!\pm\! 2\hbar n,-\!1\!\mp\! 2\hbar n)
\over B(\pm 2\hbar ,\mp 2\hbar)}
{\scriptstyle F\left(1\pm 2(n\!+\!1)\hbar,1\mp 2\hbar(n\!-\!1),2\pm 
2\hbar n;x\right)}\right]       
\stackrel{(\ref{value}),(\ref{betagamma})}{\approx}
x^{\mp 2\hbar}
{\Gamma(1\!\pm\! 2\hbar n)\,\Gamma(1\!\pm\! 2\hbar n)\over
\Gamma(1\!\pm\! 2\hbar (n\!+\!1))\,\Gamma( 
1\!\pm\! 2\hbar (n\!-\!1))}
\end{array}
\]
\[
\begin{array}{lcl}
f_1(x)
&\stackrel{(\ref{deri})}{=}&
\frac{2\hbar}{1\pm 2\hbar n}F(1\!\pm\! 2\hbar,1\! \mp\!  
2\hbar,2\!\pm\! 2\hbar n;1\!-\!x) 
x^{\mp 2\hbar}(1-x)^{1\pm2\hbar(n\!-\!1)}\nonumber\cr
&\stackrel{(\ref{lintra})}{=}&
\frac{2\hbar}{1\!\pm\! 2\hbar n}
x^{\mp 2\hbar} (1\!-\!x)^{1\!\pm\!2\hbar(n\!-\!1)}\!
\left[\!{B(2\!\pm\! 2\hbar n,\!\pm\! 2\hbar n)\over
B(1\!\pm\! 2\hbar (n\!+\!1),1\!\pm\! 2\hbar (n\!-\!1))}
{\scriptstyle  F(1\!\pm\! 2\hbar,1\!\mp\! 2\hbar,1\!\mp\! 
2\hbar n;x)}\right. \nonumber\cr
& & + \left.
x^{\pm 2\hbar n}{B(2\!\pm\! 2\hbar n,\mp\! 2\hbar n)
\over B(1\pm 2\hbar ,1\mp 2\hbar)}
{\scriptstyle 
F(1\pm 2(n\!+\!1)\hbar,1\mp 2\hbar(n\!-\!1),1\pm 2\hbar n;x)}
\right]        \nonumber\cr
&\stackrel{(\ref{value})}{\approx}& 
\frac{2\hbar}{1\pm 2\hbar n}
x^{\mp 2\hbar} \left[{B(2\!\pm\! 2\hbar n,\!\pm\! 2\hbar n)\over
B(1\!\pm\! 2\hbar (n\!+\!1),1\!\pm\! 2\hbar (n\!-\!1))} +
x^{\pm 2\hbar n}{B(2\!\pm\! 2\hbar n,\mp\! 2\hbar n)
\over B(1\pm 2\hbar ,1\mp 2\hbar)}\right]
\nonumber\cr
&\stackrel{(\ref{betagamma})}{=}&
 x^{\mp 2\hbar}
\left[\pm\frac 1n{\Gamma(1\!\pm\! 2\hbar n)\Gamma(1\!\pm\!
 2\hbar n)\over\Gamma(1\!\pm\! 2\hbar (n\!+\!1))\Gamma(
1\!\pm\! 2\hbar (n\!-\!1))} +2\hbar x^{\pm 2\hbar n}
{ \Gamma(1\!\pm\! 2\hbar n)\Gamma(\mp\! 2\hbar n)
\over  \Gamma(1\pm 2\hbar)\Gamma(  1\mp 2\hbar)}\right].
\end{array}
\nonumber
\]
For the combination $nf_1\!\mp\!f_2$ we thus find
\eqa
& nf_1\!\mp\!f_2 &\stackrel{x\rightarrow 0}{\approx}
 2\hbar x^{\pm 2\hbar (n\!-\!1)}
{ \Gamma(1\!\pm\! 2\hbar n)\Gamma(\mp\! 2\hbar n)
\over  \Gamma(1\pm 2\hbar)\Gamma(  1\mp 2\hbar)}
   \stackrel{(\ref{betagamma})}{=}\mp n  
x^{\pm 2\hbar (n\!-\!1)}
{ \Gamma(   \pm\! 2\hbar n)\Gamma(\mp\! 2\hbar n)
\over  \Gamma( \pm 2\hbar)\Gamma(   \mp 2\hbar)}
\nonumber\cr
& &    \stackrel{(\ref{bingo})}{=}\mp
\frac{q-q^{-1}}{q^n-q^{-n}}\,x^{\pm 2\hbar (n\!-\!1)}
= \mp{1\over [n]_q}\,x^{\pm 2\hbar (n\!-\!1)}.
\ena

The limits $l_i$ are thus given by
\eqa
l_1 & =&  \mp{1\over [n]_q}\nonumber \\
l_3 & =& \pm \, 1 \\
l_2 & =& {n\over n\!+\! 1}l_3-{1\over n\!+\! 1}l_1
=\pm {n\over n\!+\! 1}\left(1+{1\over [n]_q}\right),
\nonumber
\ena
which plugged into eq. (\ref{brodo}) give
\eqa
wa^+\cdot(M^{-1}\,V\, M)\,va  
&=& q^{\mp n} \left[-\frac{(w\!\cdot\!v)}{[n]_q}\,
+(w\!\cdot\!a)(a^+\!\!\cdot\! v)\frac n{n+1}\left(
\frac 1{[n]_q} +  q^{\pm (n\!+\!1)}\right)\right]
\nonumber \\
&=& \mp {q^{\mp 1}\over n_{q^{\pm 2}}}
\pm (w\!\cdot\!a)(a^+\!\!\cdot\! v)
q^{\mp 1}{(n\!+\!1)_{q^{\pm 2}}\over (n\!+\!1)}
{n\over n_{q^{\pm 2}}};
\label{brodino}
\ena
eq. (\ref{figat3}) is manifestly satisfied once we
replace the latter result in it. $\Box$

\subsection*{Acknowledgments}
It is a pleasure to thank J.~Wess  for his 
stimulating remarks  in the early 
stages of this work, as well for the 
warm hospitality at his institute.
I also thank Prof. B. Zumino
for his warm hospitality during a short visit at LBNL,
where some involved computation of this work was
completed.
This work was supported in part through a TMR fellowship
granted by the European Commission, Dir. Gen. XII for Science,
Research and Development, under the contract ERBFMICT960921, 
and in part through the U.S. Dep. of Energy under contract
No. DE-AC03-76SF00098 (during my visit at LBNL).

\end{document}